\documentclass{article}

\usepackage{amssymb, amsthm} 

\newtheorem{prop}{Proposition}[section]
\newtheorem{thm}[prop]{Theorem}
\newtheorem{cor}[prop]{Corollary}
\newtheorem{lem}[prop]{Lemma}
\newtheorem{conj}[prop]{Conjecture}
\newtheorem{algo}[prop]{Algorithm}

\newcommand{\card}[1]{{\arrowvert #1 \arrowvert\,}}
\newcommand{\class}[1]{{\mathbb{Q}_{\, class} #1}}
\newcommand{\set}[1]{{\{ #1\}}}
\newcommand{\st}{{\,\arrowvert\,}}
\newcommand{\ga}[1]{{\mathbb{Q} #1}}

\newcommand{\ind}[1]{{ 1_{W_{#1}}^{W} }}

\begin{document}
\title{Properties of the Solomon algebra homomorphism}
\author{{\bf Christophe Hohlweg}\\ 
   Institut de Recherche Math\'ematique Avanc\'ee,\\
   Universit\'e Louis Pasteur et CNRS,\\
   7 rue Ren\'e Descartes,\\
   67084 Strasbourg Cedex, France\\
   E-mail: hohlweg@math.u-strasbg.fr}

\maketitle

\tableofcontents

\section{Introduction}

Let $(W,S)$ be a finite Coxeter system. Let $I\subseteq S$ and $W_I$ be
the \textit{parabolic subgroup} of $W$ generated by $I$. We denote $X_I$ the cross section of
$W/W_I$ consisting of the unique coset representatives of minimal length.
Its elements are called \textit{the minimal coset representatives} (see \cite{bourbaki}, \cite{humphreys}, 5.12).

 The \textit{Solomon descent algebra} $\Sigma W$ (first defined in \cite{solomon}) is the subalgebra of $\mathbb Q W$ generated by $\{ x_I\}_{I\subset S}$,
 where $x_I = \sum_{w\in X_I} w$. 

Solomon has also  shown  that 
the linear function $\Phi$ which maps $x_I$ to $\ind I$, the induced character 
of the trivial representation on $W_I$ to $W$, is an algebra homomorphism called \textit{ the Solomon homomorphism}.\\

For the symmetric groups (which are Coxeter groups of type $A$), properties of $\Phi$ are related to enumerative results: certain joint statistics on the symmetric group may be enumerated by the scalar products of appropriate characters.

For example, Gessel  has shown (\cite{gessel}) that the number of element in $X_I \cap X^{-1}_J$ is the scalar product of the characters $\ind I$ and $\ind J$. We show that this result is equivalent to the following:   $\Phi$ is an isometry for a particular scalar product on $\ga S_n$ (defined by $<w,v>=1$ if $w=v^{-1}$, $=0$ otherwise).  In Section~\ref{idem},  we  extend this result of Gessel on symmetric groups to all finite Coxeter groups.

 As another example, there is a set of orthogonal idempotents $E_\lambda \in \Sigma S_n$, indexed by the  partitions $\lambda$ of $n$, which decompose $\ga S_n$ (see \cite{garsia}). Gessel and Reutenauer have shown (\cite{gesselreut}) that the scalar product of $\ind I$ and $\chi_\lambda$ - the character of the left action of $S_n$ on $\ga S_n E_\lambda$ - is the number of elements in $X_I$ of cycle type $\lambda$ (equivalently, in the conjugation class corresponding to $\lambda$. By this way, J\" ollenbeck and Reutenauer in 
\cite{jollen} have shown that $ \Phi$ has the following symmetry property:
 $$(\star)\qquad\Phi(x)(y) = \Phi(y)(x)$$ 
for all $x,y\in \Sigma S_n$ ($\Phi(x)$ is linearly extended to $\ga S_n$). In other words
$$
\sum_{w\in X_I} \ind J (w) = \sum_{g\in X_J} \ind I (g)
, \textrm{ for all } I,J\subset S .$$
Conversely, the symmetry property implies the previous enumerative result.\\

In our attempt to extend to all  finite Coxeter groups the symmetry property $(\star)$, we use the work of Bergeron, Bergeron, Howlett and Taylor (\cite{bergeron}): they have found for each finite Coxeter group a set of orthogonal idempotents $E_\lambda$, indexed by Coxeter classes  which decompose $\ga W$ (a Coxeter class is an equivalence class of subsets of $S$ conjugated under $W$).  In  Section~\ref{sectionconj}, we prove the scalar product of $\ind I$ and $\chi_\lambda$ - the character of the right action of $W$ on $E_\lambda \ga W$ - is equal to the number of  elements in $X_I$ of Coxeter type $\lambda$ if and only if $\Phi$ has the symmetry property $(\star)$ (the Coxeter type may be viewed as the  cycle type of a permutation). Thus, again, properties of $\Phi$ are strongly related to enumerative results.

We show that $\Phi$ has  property $(\star)$  for all finite Coxeter groups whose irreducible components are of type  $A_n$, $E_6$, $E_7$, $E_8$, $F_4$, $H_3$, $H_4$ or $I_2 (m)$. Whether this is true for all finite Coxeter groups is an open conjecture. Moreover, we state a stronger conjecture 
 (obtained jointly with M. Schocker) that depends on two subsets $I,J$ of $S$; we prove this new conjecture in the following cases: $I = \emptyset$, $I= S$ and $I$ is a singleton.  In particular, the new conjecture is true for dihedral groups. The general case remains open.

As a byproduct of our techniques, in Section~\ref{sectionhom},  we give a quick proof of the following result: For each conjugacy class $C$, 
if $I, J$ are conjugated under $W$, then $\card{C\cap X_I} = \card{C\cap X_J}$. 
This result has been obtained independently for all Coxeter groups by Fleischmann in \cite{fleisch}. We introduce the image of $\Phi$
 in the algebra of class functions and we describe it using some results of Bergeron, Bergeron, Howlett and Taylor (\cite{bergeron}). Moreover, we show that $\Phi$ is surjective if and only if $W$ is a direct product of symmetric groups.

\section{The Solomon  homomorphism}\label{sectionhom}

\subsection{The Solomon descent algebra}\label{secalg}

Let $w\in W$. We denote by $\ell(w)$ 
the length of $w$ as a word in the elements of $S$. The \textit{ascent set of} $w$ is the set
$$
S(w)=\set{s\in S\st \ell(ws)>\ell(w)}\ .
$$
The \textit{descent set of} $w$ is the set
$$
D(w)=\set{s\in S\st \ell(ws)<\ell(w)}\ .
$$
Observe that $D(w)\amalg S(w) =S$. Let $I\subset S$, then the set of minimal coset representatives $X_I$   can also be viewed as
\begin{eqnarray*}
X_I&=& \{x\in W\ \arrowvert\ \ell(xs)>\ell(x),\ \forall s\in I\}\\
&=&\set{w\in W\st S(w)\supset I}\\
&=&\set{w\in W\st D(w)\subset S-I}\ .
\end{eqnarray*}
\textbf{Parabolic components:} Any $w\in W$ has a unique factorization 
$
w=w^Iw_I
$
such that $w^I\in X_I$ and $w_I\in W_I$. Moreover,
$\ell(w)=\ell(w^I)+\ell(w_I)$ (\cite{humphreys}, 5.12). We call the couple $(w^I,w_I)$ the \textit{parabolic components} of $w$.\\

Let $I,J\subseteq S$ and $X_I^{-1}=\{x^{-1}\,\arrowvert\, x\in X_I\}$.  
Then the set
$$
X_{IJ}=X_I^{-1}\cap X_J
$$
is a cross section of the double cosets $W_{I}\backslash W /W_{J}$,
 consisting again of the unique representatives of minimal length (\cite{solomon}, Lemma $1$).

Following \cite{bergeron}, Section $2$, for $I,K\subset S$, we denote $X_{I\cap K}^{I} = W_I \cap X_K$. Then $X_{I\cap K}^{I}$
 is a cross section of $W_I /W_{I\cap K}$ 
which is constituted by the minimal coset representatives of $W_{K\cap I}$ in $W_I$. Finally, we have the following properties:
\begin{enumerate}
\item for $K\subset I\subset S$, $X_I X^{I}_K = X_K$;
\item for $I,J\subset S$, $X_J$ is the disjoint union of $X^{I}_{I\cap bJb^{-1}}b$ for all $b\in X_{IJ}$ and 
$W_I \cap b W_J b^{-1} = W_{I \cap b J b^{-1}}$.
 In particular, for all $w\in W$, there is
 a unique $b\in X_{IJ}$ such that $w^{J} \in X^{I}_{I\cap bJb^{-1}}b$;
 \item for $K,L\subset S$ and $b\in X_{KL}$ such that $bLb^{-1}=K$, we have $X_K b = X_L$.
\end{enumerate}
\textbf{Double parabolic components:} We have an analogue decompostion in parabolic components for 
double cosets: any $w\in W$ has a unique factorization
$
w=abc
$
such that $b\in X_{IJ}$, $a\in X_{I\cap bJb^{-1}}^I$ and $c\in W_J$. Moreover, $\ell(w)=\ell(a)+\ell(b)+\ell(c)$ (\cite{solomon}, Lemma $1$). 
We call the triplet $(a,b,c)$ the \textit{double parabolic components} of $w$.\\

 For $I\subset S$, 
we consider the element $x_I=\sum_{w\in X_I} w$ in $\ga W$. Solomon has shown in \cite{solomon} that the 
vector subspace $\Sigma W$  of $\ga W$ generated by these elements is a subalgebra. Moreover $\{x_I\}_{I\subset S}$ 
is a basis that satisfies
$$
x_Ix_J=\sum_{K\subset S} a_{IJK}x_K ,
$$
where $a_{IJK}=\arrowvert\{w\in X_{IJ}\,\arrowvert\, w^{-1}W_I w\cap W_J = W_K\}\arrowvert$.

Following the case of  $A_n$ and $B_n$ (symmetric and hyperoctahedral groups, 
see \cite{reutenau},  \cite{nantel}), we call $\Sigma W$ \textit{the Solomon descent algebra}.
Note that in these two case, the basis 
$
\tilde{x}_I = x_{S-I}= \sum_{D(w)\subset I} w 
$ 
is more commonly used. Note also that $\dim \Sigma W=2^{\card S}$.\\

We use a  basis of $\Sigma W$ which is useful to enumerate the elements of $W$ which have a 
given descent or ascent set. Let $I,J\subset S$. The set of elements of $W$, 
which have as ascent set  
$I$, is $Y_I=\set{w\in W\st S(w)=I}$.
We consider the following element in $\ga W$:
$$
y_I=\sum_{w\in Y_I} w\ .
$$
Then
$$
x_J = \sum_{I\supset J} y_I
$$
and, by inclusion-exclusion,
$$
y_I=\sum_{J\supset I} (-1)^{\card{I-J}} x_I\ .
$$
Therefore $\set{y_I}_{I\subset S}$ is a basis of $\Sigma W$.

\subsection{Coxeter classes}

Let $I,J\subset S$. We say that $I$ and $J$ \textit{are $W$-conjugate} 
if there is $w\in W$ such that $wIw^{-1}=J$. In this case, we write $I\sim_W J$.

The relation $\sim_W$, on the set of subsets of $S$, is an equivalence relation. For this relation, an equivalence class  
is called a \textit{Coxeter class}. The set of Coxeter classes of $W$ is denoted $\Lambda(W)$. For $I\subset S$, write $\lambda(I)$ its 
Coxeter class.

An element $c\in W$ with $c=s_1\dots s_n$ ($s_i \in S$) is a \textit{Coxeter element} if the $s_i$ are all distinct 
and if $S=\set{s_1,\dots ,s_n}$ ($c$ is the product of all elements of $S$, without repetition). We write $c_I$ 
for a Coxeter element of $W_I$. It is well-known that the following statements  give three equivalent characterisations of the $W$-conjugation (
See
 \cite{bergeron}, Proposition $4.4$ and \cite{geck}, Proposition $3.1.15$):
\begin{enumerate}
\item $I\sim_W J$,
\item the subgroups $W_I$ and $W_J$ are conjugated,
\item the Coxeter elements $c_I$ of $W_I$ and $c_J$ of $W_J$ are conjugated.
\end{enumerate}

Denote by $C(w)$ the conjugacy class of an element
 $w\in W$ and $Cl(W)$ the set of conjugacy classes of $W$.
A \textit{cuspidal class} is a conjugacy class $C$ such that $C\cap W_I=\emptyset$ for all
 $I\subsetneq S$. A \textit{$I$-cuspidal class} is a cuspidal class of the parabolic subgroup 
$W_I$. As example, the Coxeter elements are all conjugated and their conjugacy class is cuspidal (see \cite{geck},
 Proposition $3.1.6$). 
 
An element $w\in W$ is an \textit{$I$-element} ($I\subset S$) if $C(w)\cap W_I$ is a $I$-cuspidal
class. If $J\in \lambda(I)$, then $w$ is also a $J$-element (since $W_J$ is conjugated
 to $W_I$, $C(w)\cap W_I$ is conjugated to $C(w)\cap W_J$. As subgroups of $W_J$ are conjugated to
 subgroups of $W_I$, $C(w)\cap W_J$ is a $J$-cuspidal class).
 
More generally, for any $\lambda\in \Lambda$, we say that $w\in W$ is \textit{a $\lambda$-element} (or is of \textit{Coxeter type} $\lambda$)
 if $w$ is a $I$-element for some $I\in \lambda$. Denote by $C(\lambda)$ the set of 
 $\lambda$-elements. We call it \textit{the $\lambda$-set}. 

It is well-known that  conjugacy classes and  Coxeter classes coincide in the symmetric group. 
The subclass of finite Coxeter groups that have this property is given by the following proposition. We will see that this 
is also the subclass of finite Coxeter groups for which the Solomon homomorphism is surjective.

\begin{prop}\label{coxeqequ} Let $(W,S)$ be a finite Coxeter system then the following 
propositions are equivalent:
\begin{enumerate}
\item $W$ is a direct product of symmetric groups;
\item $\card{\Lambda(W)}=\card{Cl(W)}$;
\item Coxeter classes coincide with conjugacy classes.
\end{enumerate}
\end{prop}

We need two lemmas to prove Proposition~\ref{coxeqequ}. The first one implies that there are more conjugacy classes than Coxeter classes.

\begin{lem}\label{coxconjclass} Let $\lambda\in \Lambda$, then
 $
  C(\lambda)= \amalg C
 $, 
where the disjoint union is taken over all $C\in Cl(W)$ such that $C\subset C(\lambda)$. 
Moreover, $W$ is the disjoint union of its $\lambda$-sets, for $\lambda\in \Lambda$.
\end{lem}
\begin{proof} If $w\in C(\lambda)$, then $C(w)\cap W_I$ is a cuspidal class of $W_I$ 
(for some $I\in \lambda$). Thus all $g\in C(w)$ are $\lambda$-elements. 
Therefore $C(w)\subset C(\lambda)$. Finally, for all $w\in C(\lambda)$, there is a conjugacy
 class $C$ such that $w\in C \subset C(\lambda)$ and the lemma is proved. 
\end{proof}

It is well-known that the conjugacy classes of $W$ are the direct product of conjugacy classes of its ireducible components (as $W$ is isomorphic to
  the direct product of them). This second lemma gives an analogue for Coxeter classes.
  
\begin{lem}\label{coxclass1} Let $(W,S)$ be a finite Coxeter system. For $i=1,\dots n$, let $(W_i , K_i)$ be the irreducible
 components of $(W,S)$. For $\lambda\in\Lambda (W)$ and $I\in \lambda$, denote $I_p = I \cap K_p$ and $\pi_p (\lambda) = \lambda (I_p)\in\Lambda(W_p)$. 
 Then the following mapping
 \begin{eqnarray*}
 \pi \, :\ \Lambda (W) &\longrightarrow & \Lambda (W_1) \times \dots \times \Lambda (W_n) \\
\lambda & \longmapsto& (\pi_1 (\lambda),\dots ,\pi_n (\lambda))
 \end{eqnarray*}
is a bijection. Moreover, for any $\lambda\in\Lambda (W)$
 $$
 \card{C(\lambda)} = \prod_{p=1}^{n}  \card{C(\pi_p (\lambda))} .
 $$
\end{lem}
\begin{proof} It is readily seen that one just has to prove the lemma in the case $n=2$ (and
 conclude by induction on $n$). As $W_p$, $p=1,2$, are the irreducible 
components of $W$, they can be viewed as normal parabolic subgroup of $W$ such that $W= W_1 W_2$. Therefore, for any $w\in W$,
 $w = w_1 w_2 = w_2 w_1$ with $w_p \in W_p$. Observe that $w^{K_1} = w_{K_2} = w_2$ and $w_{K_1} = w^{K_2} = w_1 $. For $I\subset S$ and 
 $p = 1,2$, $W_p \cap I \subset W_p \cap S = K_p$. 
Thus $I_p = I\cap K_p = I\cap W_p$ and $I$ is the disjoint union of $I_1$ and $I_2$.
 Finally, since each element of $W_1$ commutes with each element of $W_2$,
 $w I_p w^{-1} = w_p I_p w_{p}^{-1}$. 
  
One first shows that $\pi$ is well defined. Let $I,J\in\lambda$, there is $w\in W$ such that $wIw^{-1} = J$.  As $W_p$ is a normal subgroup, for $p=1,2$,
$w I_p w^{-1} =J_p$. Moreover, since $w I_p w^{-1} = w_p I_p w_{p}^{-1}$, one has $w_p \in W_p$ such that $ w_p I_p w_{p}^{-1} =J_p $,  
hence $\lambda (I_p) = \lambda (J_p) \in \Lambda (W_p)$. Then  $\pi_p$, and therefore $\pi$, are well defined. 

Let $(\lambda_1 ,\lambda_2)\in \Lambda (W_1) \times \Lambda (W_2)$ and $I_p \in \lambda_p$. Then $I = I_1 \cup I_2 \subset S$ and 
$\lambda (I) \in\Lambda (W)$. It is readily seen that $\pi_p (\lambda(I)) = \lambda_p$ hence $\pi (\lambda(I) =(\lambda_1 ,\lambda_2)$.
 Therefore, $\pi$ is surjective.
 
 Let $\lambda ,\mu \in \Lambda (W)$ such that $\pi (\lambda) = \pi(\mu)$. So $\pi_p (\lambda) = \pi_p (\mu)$. Let $I\in \lambda$ and $J\in \mu$, 
 then $I_p$ and $J_p$ are $W_p $-conjugate by an element $w_p \in W_p$. Write $w = w_1 w_2 = w_2 w_1$, then one has $w I w^{-1}= J$.
Therefore, $I\sim_W J$ so $\lambda = \mu$. One concludes that $\pi$ is injective, hence bijective.

Let $\lambda\in\Lambda (W)$ and $I\in\lambda$. If $w\in C(\lambda)$, then $C(w)\cap W_I$ is a $I$-cuspidal class. it is well-known that
$C^{p}(w_p)$ is a conjugacy class of $W_p$ and $C^{p}(w_p) \cap W_{I_p}$ is an $I_p $-cuspidal class (see \cite{geck}, exercise $3.10$). Therefore, 
the mapping
 \begin{eqnarray*}
 \eta \, :\ C(\lambda) &\longrightarrow & C(\pi_1 (\lambda)) \times C(\pi_2 (\lambda))\\
 w & \longmapsto& (w_1 , w_2)
 \end{eqnarray*}
is well defined and injective since the parabolic components $(w_1 , w_2)$ of $w$ are uniquely determinated. Moreover, $\eta$ is easily seen to be surjective, hence
the lemma is proved.
\end{proof}

\begin{proof}[Proof of Proposition~\ref{coxeqequ}] One  may assume that $W$ is irreducible, by Lemma~\ref{coxclass1}.
\begin{enumerate}
 \item $\card{\Lambda(W)}\leq\card{Cl(W)}$, and $\card{\Lambda(W)}<\card{Cl(W)}$ if $W$
 contains a parabolic subgroup $W_I$ that has at least two cuspidal classes. Indeed, denote
 by $C_1^I$ and $C_2^I$ two distinct cuspidal classes of $W_I$. Let $w_1\in C_1^I$ and 
 $w_2\in C_2^I$.  Recall a useful result about cuspidal classes:
 if $I\subset S$, if
 $C^I \in Cl(W_I)$ is a $I$-cuspidal class and if $w\in C^I$, then
$C^I=C(w) \cap W_I$ (see \cite{geck}, Theorem $3.2.11$). This implies that the conjugacy classes $C(w_1)$ and $C(w_2)$ of $W$ are distinct
 and contained in $C(\lambda(I))$. Therefore $\card{\Lambda(W)}<\card{Cl(W)}$.
 \item The Coxeter group of type  $I_2(m)$, $m\geq4$ has at least two cuspidal
 classes. Indeed, let $S=\set{s,t}$, then the conjugacy class $C_1$ of the Coxeter element 
$st$ is cuspidal. Moreover,
  $stst$ is an element of minimal length in its conjugacy class $C_2$. 
As $\ell (st) =2 < 4=\ell(stst)$, $C_2$ is distinct
 of $C_1$. Observe that $C_2$ is cuspidal, hence $C_1 , C_2$ are two distinct cuspidal classes.
 \item The Coxeter group of type $D_4$ has at least two cuspidal
 classes. Indeed, $S=\set{1,2,3,4}$ with $m(1,2)=m(2,3)=m(2,4)=3$. 
 The conjugacy class $C_1$ of the Coxeter element $1234$ is cuspidal. Moreover
  $12342$ is an element of minimal length in its conjugacy class $C_2$ thus $C_2$ is distinct
 of $C_1$. Observe that $C_2$ is cuspidal.
 \item The only irreducible Coxeter group that has no parabolic subgroup of type 
$D_4$ and $I_2(m)$ is the symmetric group (type $A$, see classification of finite 
Coxeter groups in \cite{humphreys}). \end{enumerate}\end{proof}

\subsection{The Solomon homomorphism}

 Denote by $\class W$ the set 
of class functions on $W$. In his article \cite{solomon}, Solomon shows that the following linear homomorphism
\begin{eqnarray*}
\Phi\ :\ \Sigma W &\longrightarrow & \class W \\
               x_I   &\longmapsto & 1_{W_I}^{W}
\end{eqnarray*}
 is an algebra homomorphism. Moreover he shows that
$$
\ker\Phi = Rad\ \Sigma W
$$
is generated by the elements $x_I-x_J$ with $I\sim_W J$.

\paragraph*{Remark:}\label{corsolomon} Note that the theorem implies that if $I \sim_W J$,
 then $\ind I = \ind J$; this may also be deduced from general results on induction of characters, 
since $W_I$, $W_J$ are conjugated in $W$.

The following proposition has been obtained independently for all Coxeter groups by Fleischmann (\cite{fleisch}).

\begin{prop}\label{propker} Let $I,J\subset S$ such that $I\sim_W J$ then
\begin{enumerate}
\item for any conjugacy class $C$, $\card{X_I \cap C} =\card{X_J \cap C}$;
\item for any Coxeter class $\lambda$, $\card{X_I \cap C(\lambda)} =\card{X_J \cap C(\lambda)}$.
\end{enumerate}
\end{prop}

\begin{lem}\label{hom2} Extend linearly class functions to $\ga W$. Let $x\in \ker\Phi$, then 
$$
f(x)=0,\qquad \forall f\in \class W\ .
$$
\end{lem}
\begin{proof} Let $f\in \class W$. As $f$ is a linear combination of characters, 
one can assumes that $f$ is a character, associated to
 some representation $\mu$ of $W$, of degree $n$.
 As $x\in \ker\Phi = Rad\ \Sigma W$, 
$x$ is a nilpotent element of $\Sigma W$. Thus, the characteristic polynomial of $\mu x$ is $X^{n}$.
 As $f(x) = tr(\mu x)$ is the coefficient of $\pm X^{n-1}$ in this polynomial, one has $f(x)=0$.
\end{proof}

\begin{proof}[Proof of Proposition~\ref{propker}] Let $f_C$ be the characteristic function of the conjugacy class $C$.
 It is well-known that $f_C \in\class W$. Then $f_C (x_I - x_J) = 0$, by Lemma~\ref{hom2}. Hence
 $$
 \card{X_I \cap C}=f_C (x_I ) =f_C (x_J) =\card{X_J \cap C} .
 $$
 The second statement is immediate from $1$ and Lemma~\ref{coxconjclass}.
\end{proof}

\subsection{Image of the Solomon homomorphism}

Let $\lambda$ be a Coxeter class. As $\ind I=\ind J$ for all $I,J\in\lambda$, 
we denote $\varphi_\lambda$ instead of $\ind I$ for all
 $I\in \lambda$. Following \cite{bergeron} (Theorem $4.5$), we consider the matrix
 $$
 A=\left(\varphi_\lambda (c_\mu)\right)_{\lambda,\mu\in \Lambda}\ ,
 $$
where $c_\mu$ is an arbitrary $\mu$-element (for example a Coxeter $\mu$-element). 
In \cite{bergeron}, the authors have shown that:
\begin{enumerate}
\item $A$ is an invertible matrix; 
\item $ \dim Rad \Sigma =\dim \ker\Phi= 2^{\arrowvert S\arrowvert} - \arrowvert \Lambda (W) \arrowvert$.
\end{enumerate}
They have also introduced the characteristic class function of $\lambda$-elements
$$
\xi_\lambda = \left\{\begin{array}{ll}
 1&\textrm{if}\ w\in C(\lambda) ;\\ 
0&\textrm{otherwise}\ .\end{array}\right. .
$$
 Following \cite{bergeron}, we denote by 
$\nu_{\lambda\,\mu}$ the coefficients of $A^{-1}$. Observe that if $w\in C(\mu)$, then $\xi_{\lambda} (w) = \xi_{\lambda} (c_{\mu})$, for a Coxeter
 $\mu$-element $c_\mu$. Then we have
$$
\xi_\lambda = , \sum_{\beta\in\Lambda} \nu_{\lambda\,\beta}\,\varphi_\beta\ 
$$
since $\sum_{\beta\in\Lambda} \nu_{\lambda\,\beta}\,\varphi_\beta (c_\mu) = \delta_{\lambda , \mu}$.\\

We denote $\Gamma W$ the image of $\Phi$ in $\class W$ (spanned by $\{\ind I\}_{I\subset S}$). 
As $\Sigma W$ is a $\mathbb Q$-algebra, $\Gamma W$ is also a $\mathbb Q$-algebra. 
The results of Bergeron, Bergeron, Howlett and Taylor can be rewritten as follow: $\{\varphi_\lambda\}_{\lambda\in\Lambda(W)}$, $(\xi_\lambda)_{\lambda\in\Lambda} $
 are two basis of $\Gamma W$ and $\dim\Gamma W=\arrowvert \Lambda (W) \arrowvert$. Therefore, the result of Bergeron, Bergeron, Howlett and Taylor about 
 the basis $(\xi_\lambda)_{\lambda\in\Lambda} $ implies the following proposition, which states that,  $\Gamma W$
 is the set of functions on $W$ stable on Coxeter classes.

\begin{prop}  A class function $f$ is in $\Gamma W$ if and only if $f$ is a $\mathbb Q$ valued function
 and is stable on $\lambda$-sets, for all $\lambda\in \Lambda$.
\end{prop}
\begin{proof} By Lemma~\ref{coxconjclass}, $\xi_\lambda=\sum_{C\subset C(\lambda)} f_C$ 
where $f_C$ is the characteristic function of the 
conjugacy class $C$. If $f$ is stable on $\lambda$-set, 
 for any $\lambda$, the coefficients of $f$ in the basis $\set{f_C}_{C\in Cl(W)}$ are equal
 if the conjugacy classes are in the same $\lambda$-set. Therefore $f$ may be written in the 
basis $\set{\xi_\lambda}_{\lambda\in\Lambda W}$ with coefficients in $\mathbb Q$ if 
$f$ is a $\mathbb{Q}$ valued function. On the other side, if $f\in\Gamma W$, it is readily seen that $f$ is $\mathbb Q $-valued and is stable on all $\lambda$-sets.
\end{proof}

We conclude  the study of the image of $\Phi$ with the following corollary.

\begin{cor}[of Proposition~\ref{coxeqequ}]\label{surj} Let $(W,S)$ be a finite Coxeter system. Then
 $\Phi$ is surjective if and only if $W$ is a direct product of symmetric groups.
\end{cor}

\section{Character formula and enumeration}\label{idem}

\subsection{The Solomon homomorphism is an isometry}

It is well known that irreducible characters of $W$ are $\mathbb R$-valued. More precisely,
 there is a finite field extension $\mathbb K$ of $\mathbb Q$ such that the characters of $W$ can be
 realized over $\mathbb K$ (see \cite{geck}, Chapter $5$). For example, for type $A_n$, $B_n$ 
and $D_n$, $\mathbb K=\mathbb Q$. As a well-known consequence, each element of $W$ is
 conjugated to its inverse and  the classical scalar product 
  on $\class W$  is given by
$$
<f,g>_W = \frac{1}{\card W} \sum_{w\in W} f(w) g(w) .
$$
We consider the following scalar product on $\mathbb Q W$. For $w,g\in W$ let
$$
<w,g> = \left\{\begin{array}{ll}
 1&\textrm{if}\ w=g^{-1} ;\\ 
0&\textrm{otherwise} .\end{array}\right.
$$

This scalar product is of great interest concerning the problems of enumeration:
 Let $M\subset W$ and $z_M = \sum_{w\in M} w \in \ga W$. Then the  number of elements of $W$ which 
are in $N$ and whose inverse 
is in $N\subset W$ is equal to $<z_M , z_N>$.

For example, let $I,J\subset S$; then
\begin{enumerate}
\item  $<x_I,x_J> = \card{X_{IJ}}$ is the number of $w$ in 
$W$ such that the ascent 
set of $w$ contains $I$ and that the ascent set of $w^{-1}$ contains $J$;
\item $<y_{I},y_{J}> $ is the number of $w$ in $W$ such that the ascent 
set of $w$ equals $I$ and that the ascent set of $w^{-1}$ equals $J$, that is,  
the number of $w$ in $W$ such that the descent 
set of $w$ equals $S - I$ and that the descent set of $w^{-1}$ equals $S-J$.
\end{enumerate}

The following theorem gives an isometry property of $\Phi$ for this scalar product and 
the classical scalar product on $\class W$. This is a generalization of Gessel's result for 
symmetric groups (see \cite{gessel}).

\begin{thm}\label{isom} Let $x,y\in \Sigma W$, then 
$<\Phi(x),\Phi(y)>_W\, =\, <x,y>$.
\end{thm}
\begin{proof} As $\{x_I\}_{I\subset S}$ is a basis of $\Sigma W$, one just has to prove 
the theorem for these elements.
Let $I,J\subset S$, let $\chi_I$ (resp.  $\chi_J$) be a character on
 $W_I$ (resp. $W_J$).  Denote  $\chi_I^W$ (resp.  $\chi^W_J$) be the character induced by
 $\chi_I$ (resp.  $\chi_J$) on $W$.  As $X_{IJ}$ is a set of representatives of double coset $W_I w W_J$, Mackey's Theorem (see \cite{ledermann}, 
Theorem $3.3$, p.$86$) implies 
$$
\chi_I^W \chi_J^W = \sum_{x\in X_{IJ}} \beta_{x}^{W} ,
$$
where $\beta_{x}^{W}$ denote the induced character on $W$ of the character $\beta_x$ on the subgroup
$W_I \cap x W_J x^{-1} = W_{I\cap xJx^{-1}}$ defined by
$$
\beta_{x} (w) = \chi_I (w) \chi_J (x^{-1}wx)  ,
$$
for $w\in W_{I \cap xJ x^{-1}}$. 
Applying this to $\chi_I = 1_{W_I}$ and $\chi_J =  1_{W_J}$, one has 
$$
1_{W_I}^{W}  1_{W_J}^{W} = \sum_{x\in X_{IJ}} 1^{W}_{W_{I \cap x J x^{-1}}}   .
$$
Therefore
\begin{eqnarray*}
<\Phi(x_I),\Phi(x_J)>_W &=& <1_{W_I}^{W} 1_{W_J}^{W}, 1>_W\\
&=&\sum_{x\in X_{IJ}} <1_{W_{I \cap x J x^{-1}}}^{W},1>_{W}\\
&=&\sum_{x\in X_{IJ}} <1,1>_{W_{I \cap x J x^{-1}}} , \ \textrm{by Frobenius reciprocity}\\
&=& \card{X_{IJ}}\\
&=&<x_I , x_J >  .
\end{eqnarray*}
\end{proof}

\begin{cor}[of proof] Let $I,J\subset S$, then
$$
\arrowvert X_{IJ} \arrowvert = <1_{W_I}^{W},1_{W_J}^{W}>_W .\\
$$
Equivalently, the number of $w$ in 
$W$ such that the ascent 
set of $w$ contains $I$ and that the ascent set of $w^{-1}$ contains $J$ is equal 
to $<1_{W_{I}}^{W},1_{W_{J}}^{W}>_W$.
\end{cor}

\subsection{Idempotents in the group algebra}\label{sec32}

Let $e=\sum_{w\in W} e_{w} w\in \mathbb Q W$ be an idempotent. Then $e\,\mathbb Q W$ (resp. 
$\mathbb Q W e$) is a right (resp. left) $W$-module by right (resp. left) multiplication. 
Denote by $\chi_ R$ (resp. $\chi_ L$) the character of the representation 
$e\,\mathbb Q W$ (resp. $\mathbb Q W e$).

\begin{prop} The right $W$-module $e\,\mathbb Q W$ and the left $W$-module $\mathbb Q W e$ have 
the same character  $\chi_e (w)=\chi_L (w) = \chi_R (w)$. In particular
 $\dim \mathbb Q W e = \dim e\,\mathbb Q W$.
\end{prop}

The proposition is an immediate corollary of the following well-known lemma. For a proof, see
 \cite{reutenau}, Lemma $8.4$, where it is given for any finite group and for 
the classical scalar product on $\ga W$ (for which $W$ is an orthonormal basis). 

\begin{lem}\label{idem1} For any idempotent $e\in \mathbb Q W$  and any $w\in W$ we have
$$
\chi_L (w) = \chi_R (w) = \sum_{x\in W} < x^{-1} e x, w> . 
$$
\end{lem}
The next well-known lemma
is a important tool for the following.
\begin{lem}\label{idem2} For any idempotents $e$ and $e'$ in $\ga W$ we have:
$$
<\chi_e , \chi_{e'}>_W = \chi_{e'} (e) .
$$
\end{lem}
\begin{proof} 
\begin{eqnarray*}
<\chi_e,\chi_{e'} >_W &= &\frac{1}{\arrowvert W\arrowvert}\sum_{w\in W} \chi_e (w)\chi_{e'}  (w)\\
&=&\frac{1}{\arrowvert W\arrowvert}\sum_{w\in W}\sum_{x\in W} <x^{-1}ex,w>\chi_{e'}  (w)\, ,\qquad \textrm{by Lemma~\ref{idem1}}\\
&=&\frac{1}{\arrowvert W\arrowvert}\sum_{x\in W} \sum_{w\in W} e_{x w^{-1} x^{-1}} \chi_{e'}  (w)\\
&=&\frac{1}{\arrowvert W\arrowvert}\sum_{x\in W} \sum_{w\in W} e_{x w^{-1} x^{-1}} \chi_{e'}  (w^{-1}), \textrm{ since } w\sim w^{-1} \\
&=&\frac{1}{\arrowvert W\arrowvert}\sum_{x\in W} \sum_{w\in W} e_{x w^{-1} x^{-1}} \chi_{e'}  (x w^{-1} x^{-1})\\
&=&\frac{1}{\arrowvert W\arrowvert}\sum_{x\in W} \sum_{g\in W} e_{g} \chi_{e'}  (g)\\
&=&\chi_{e'}  (e) .
\end{eqnarray*}
\end{proof}

For $I\subset S$, $e_{I}= \frac{1}{\card{W_I}} \sum_{w\in W_I} w$ is an idempotent of $\mathbb Q W$. Observe that
$\chi_{e_I}=\ind I$. Let 
$$
Fix_{W_I} (e\ga W) = \set{x\in e\ga W \st xu = x\, ,\ \forall u\in W_I} .
$$
Observe that $e\,\mathbb Q W$ is a right $W_I$-module and that $W_I$ acts trivially
on $e\,\mathbb Q W e_I$. Therefore 
$$
Fix_{W_I} (e\ga W) = e\,\mathbb Q W e_I 
$$
and 
$$
 \dim Fix_{W_I} (e\ga W)= <\chi_e , \ind I>_W  ,
$$
see \cite{fulton}, p.$15-16$. Hence 
$$
\dim e\,\mathbb Q W e_I = <\chi_e , \ind I>_W =\chi_e (e_I)  .
$$

In \cite{solomon68}, the author has looked at the following idempotents: Let $I\subset S$ and $\epsilon$ be the sign character, denote
$$
\tilde{e}_I = \frac{1}{\card{W_I}}\sum_{w\in W_I} \epsilon (w) w.
$$
Solomon has shown in \cite{solomon68}  that the character $\tilde{\chi}_I$ of the action of $W$ on $\tilde{e}_{S-I}e_I \ga W $ is $\Phi(y_I)$. 
We deduce, from Theorem~\ref{isom}, the following corollary which generalize to finite Coxeter groups a Gessel's
 result for symmetric groups (\cite{gessel}) and a Poirier's result for hyperoctaedral groups 
(\cite{poirier}).

\begin{cor} Let $I,J\subset S$, then
\begin{enumerate}
\item  the number of $w$ in $X_J$ such that the descent 
set of $w^{-1}$ equals $S - I$ 
is $<\tilde{\chi}_{I} , \ind J>_W = \dim Fix_{W_J}(\tilde{e}_{S-I} e_I \ga W)$;
\item the number of $w$ in $W$ such that the descent 
set of $w$ equals $S - I$ and that the descent set of $w^{-1}$ equals $S-J$ 
is $<\tilde{\chi}_{I} , \tilde{\chi}_{J}>_W $.
\end{enumerate}
\end{cor}
\begin{proof}  
\begin{eqnarray*}
 <\tilde{\chi}_{I} , \tilde{\chi}_{J}>_W
&=& <\Phi(y_I),\Phi(y_J)>_W\\
&=&<y_I,y_J>\, , \qquad\textrm{by Theorem~\ref{isom}} .
\end{eqnarray*}
 As $<y_I,y_J>$ is the number of $w$ in $W$ such that the descent 
set of $w$ equals $S - I$ and that the descent set of $w^{-1}$ equals $S-J$, the second statement is proved. As $\Phi$ is an algebra homomorphism, the first one follows from
$$
x_J = \sum_{K\supset J} y_K .
$$ 
\end{proof}

\section{Conjectures and results}\label{sectionconj}

\subsection{The symmetry conjecture}

First, as for each $x\in \Sigma W$, $\Phi(x)$ is a class function on $W$,  we can linearly extend it  to $\mathbb Q W$.

In \cite{jollen},  J\" ollenbeck and Reutenauer prove that if $W$ is a Coxeter group of type $A_n$, 
then $\Phi(x)(y) = \Phi(y)(x)$, for all $x,y\in \Sigma W$. They sketch the following conjecture.

\begin{conj}[Symmetry conjecture]\label{conj1} Let $W$ be a finite Coxeter group. Then for all $x,y\in \Sigma W$ we have
$$
\Phi(x)(y) = \Phi(y)(x) .
$$
\end{conj}

The purpose of this section is to show the following theorem, which solves 
this conjecture for a subclass of finite Coxeter groups.

\begin{thm}\label{conjtrue} Let $W$ be a finite Coxeter group whose  irreducible components are of type
 $A_n$, $E_6$, $E_7$, $E_8$, $F_4$, $H_3$, $H_4$ or $I_2 (m)$.
Then for all $x,y\in \Sigma W$ we have
$$ 
 \Phi(x)(y) = \Phi(y)(x) .
$$
\end{thm}

\begin{cor}\label{hehe} With the same hypothesis than Theorem~\ref{conjtrue} one has for any $I,J\subset S$:
$$
\sum_{w\in X_I} 1_{W_J}^W (w) = \sum_{g\in X_J} 1_{W_I}^W (g).
$$
\end{cor}

\paragraph*{Restriction  to a basis.} In the next proposition, one readily sees that our investigation can be restricted to $x_I , x_J \in \Sigma W$, 
for $I,J\subset S$.
 It also proves Corollary~\ref{hehe}.

\begin{prop}\label{propsym} Let $W$ be a finite Coxeter group. Then the following propositions
are equivalent:
\begin{enumerate}
\item $\Phi(x)(y)=\Phi(y)(x)$, $\forall x,y\in \Sigma W $;
\item $\Phi(x_I)(x_J)=\Phi(x_J)(x_I)$, $\forall I,J\subset S$;
\item $\sum_{w\in X_I} 1_{W_J}^W (w) = \sum_{g\in X_J} 1_{W_I}^W (g)$, $\forall I,J\subset S$.
\end{enumerate}
\end{prop}

\paragraph*{Restriction to irreducible Coxeter systems.} Now, we show that it is enough to look  at irreducible finite Coxeter groups to prove the conjecture. 
Let $(W, S)$ be a finite Coxeter system with $(W_1 , I_1 ) , \dots , (W_n , I_n )$ its irreducible components. Thus $W$ is isomorphic to 
$W_1 \times \dots\times W_n$.

It is well-known  that if $\chi_i \in \class W_i$ and $\chi_j \in \class W_j$ then the tensor product 
$$
\chi_i \otimes \chi_j (w_1 , w_2 ) = \chi_i (w_1 ) \chi_j (w_j ),\ \forall (w_i , w_j )\in W_i \times W_j  .
$$
is a class function in $\class {(W_i \times W_j )}$. For any $w\in W$, there is a unique $w_i \in W_i$ such
 that $w= w_1 \dots w_n$. Then
$$
\chi_1 \otimes \dots \otimes \chi_n (w)= \chi_1 (w_1 )  \dots  \chi_n (w_n ) .
$$

Denote by $\Phi_p$ the Solomon algebra homomorphism from $\Sigma W_p$ to $\class W_p$. 

\begin{prop}\label{irred} If for $p=1 ,\dots , n$,
$$
\Phi_p (x_p)(y_p)=\Phi_p (y_p)(x_p)\ , \ \forall x_p ,y_p \in \Sigma W_p
$$
then 
$$
\Phi(x)(y)=\Phi(y)(x)\ , \ \forall x,y\in \Sigma W  .
$$
\end{prop}
\begin{proof} 
Proposition~\ref{propsym} implies that one just has to show $\Phi(x_K )(x_L )=\Phi(x_L )(x_K )$, for all $K,L \in S$.
Following  the proof of Lemma~\ref{coxclass1}, one just shows the proposition for $p=1,2$ and for $M\subset S$, one denotes
 $M_p = M \cap I_p \subset I_p$, $p= 1,2$. Remember that $X^{I_{p}}_{M_{p}} = W_{p} \cap X_{p}$.
By \cite{bergeron}, Lemma $3.1$, one has 
$$
X_M = X^{I_1}_{M_1}  X^{I_2}_{M_2}
$$
hence in $\ga W$,
$$
x_M = x^{I_1}_{M_1}  x^{I_2}_{M_2}  .
$$
It is well-known that
$$
1_{W_M}^{W} = 1_{W_{M_1 }}^{W_1}\otimes 1_{W_{M_2 }}^{W_2}  ,
$$
therefore, for $K,L\subset S$ one has
\begin{eqnarray*}
\Phi(x_K )(x_L )&=&\sum_{w\in X_L} 1_{W_K}^W (w)\\
 &=&\sum_{(w_1 , w_2) \in X^{I_1 }_{L_1} \times X^{I_2 }_{L_2}} 1_{W_{K_1}}^{W_1}\otimes 1_{W_{K_2}}^{W_2} (w_1  w_2)\\
 &=& \sum_{(w_1  , w_2) \in X^{I_1}_{L_1} \times X^{I_2 }_{L_2}} 1_{W_{K_1 }}^{W_1} (w_1 ) 1_{W_{K_2 }}^{W_2} (w_2)\\
 &=& \sum_{w_1  \in X^{I_1}_{L_1}}1_{W_{K_1 }}^{W_1} (w_1 ) \sum_{w_2  \in X^{I_2}_{L_2}}  1_{W_{K_2 }}^{W_2} (w_2)\\
 &=&  1_{W_{K_1 }}^{W_1} (\sum_{w_1 \in X^{I_1 }_{L_1}} w_1 ) 1_{W_{K_2 }}^{W_2} (\sum_{w_2 \in X^{I_2 }_{L_2}} w_2 )\\
 &=&   1_{W_{K_1 }}^{W_1} (x^{I_1}_{L_1}) 1_{W_{K_2 }}^{W_2} (x^{I_2}_{L_2})\\
 &=&  1_{W_{L_1 }}^{W_1} (x^{I_1}_{K_1}) 1_{W_{L_2 }}^{W_2} (x^{I_2}_{K_2}), \ \textrm{as}\ x^{I_p}_{K_p},x^{I_p }_{L_p }\in \Sigma W_p \\
 &=& \Phi(x_L )(x_K ) .
\end{eqnarray*}
\end{proof}

\paragraph*{Restriction to Coxeter classes.} 
\begin{prop}\label{isimj} Let $I,J\subset S$. If $I\sim_W J$ then $\Phi(x_I)(x_J)=\Phi(x_J)(x_I)$.
\end{prop}
\begin{proof} As $x_I - x_J \in\ker \Phi$, one has $\Phi(x_I) = \Phi (x_J)$. By Lemma~\ref{hom2}, one has
$\Phi(x_J)(x_J - x_I) = 0$. Finally,
$$
\Phi(x_I)(x_J) =  \Phi(x_J)(x_J) =\Phi(x_J)(x_I)  .
$$
\end{proof}

\subsection{Othogonal idempotents and character formulas}

In \cite{bergeron}, the authors have shown that for any finite Coxeter group $W$, 
there is a family of  idempotents $E_{\lambda}\in\Sigma W$ indexed by the Coxeter classes 
$\lambda\in \Lambda(W)$ satisfying $\Phi(E_{\lambda})=\xi_{\lambda}$. We call a family of \textit{good idempotents} such a family. 
We denote $\chi_{\lambda}$ the character of the right action 
of $W$ on $E_{\lambda}\mathbb Q W$. 

 Recall that if $e,e'$ are idempotents 
in $\ga W$, with associated characters $\chi$, $\chi '$, then $\dim (e\ga W e') = <\chi,\chi '>_W$. In particular, 
if $e' =e_I = \frac{1}{\card{W_I}} \sum_{w\in W_I} w$, then $\dim (e\ga W e_I ) = \dim (Fix_{W_I } (e\ga W)) = <\chi,\ind I >_W$. The purpose of this section is to show the following 
theorem and corollary. The theorem states that the symmetry 
conjecture (Conjecture~\ref{conj1}) is equivalent 
to the calculation of the cardinality of $C(\lambda) \cap X_I$ as the scalar product of $\ind I$ with $\chi_\lambda$.

\begin{thm}\label{mainidem} Let $W$ be a finite Coxeter group. 
Then the following items
are equivalent:
\begin{enumerate}
\item The symmetry conjecture holds;
\item For any family of good idempotents $(E_{\lambda})_{\lambda\in\Lambda (W)} \subset\Sigma W$,
 we have
\begin{eqnarray*}
\card{ C(\lambda)\cap X_I } &=& <\chi_{\lambda},1_{W_I}^W>_W   ,
  \end{eqnarray*}
  for all $I \subset S$ and for all $\lambda\in \Lambda(W)$.
\end{enumerate}
\end{thm}
\paragraph{Remark.} In the case $A_{n-1} = S_n$,  Coxeter classes and  conjugacy classes are in bijection. Moreover, they are indexed
by the partitions of $n$. 

In \cite{garsia} (see also \cite{reutenau}, Section $9.2$), the authors have shown the existence of  a good family in the case $A_{n-1}$. 
They have also shown that each indecomposable projective right module generated by $E_\lambda$, $\lambda$ partition of $n$, is the 
multilinear part of the subspaces generated by symmetrized Lie polynomiales of type $\lambda$.

Finally, in \cite{jollen}, the authors 
found an analog of 
Theorem~\ref{mainidem} for the case $A_{n-1}$ and they proved Theorem~\ref{conjtrue} in the case of symmetric group  using \cite{reutenau}, Corollary $9.44$
 that is they proved that  the second part of Theorem~\ref{mainidem} appear if $W=A_{n-1}$. 
This is how J\" ollenbeck and Reutenauer proved Theorem~\ref{conjtrue} in this case.\\

The following corollary, which generalizes a  result by Gessel and Reutenauer on the symmetric 
groups (\cite{gesselreut}) and Poirier for hyperoctahedral groups (\cite{poirier}), is an immediate consequence of Theorem~\ref{mainidem} 
and Theorem~\ref{conjtrue}.

\begin{cor} Let $W$ be a finite Coxeter group such that its irreducible components are of type $A_n$, $E_6$, $E_7$, $E_8$, $F_4$, $H_3$, $H_4$ or $I_2 (m)$.
Then the number of  $\lambda$-elements of $W$ which have $I$ as  ascent set is $<\chi_\lambda , \tilde{\chi}_I>_W$, where $\tilde{\chi}_I = \Phi(y_I)$.
\end{cor}

We need the following lemma to prove Theorem~\ref{mainidem}. It is the consequence of the results on $\ker\Phi$ (\cite{solomon}) 
and $Im \Phi$ (\cite{bergeron}).

\begin{lem}\label{hom1} Let $(E_{\lambda})_{\lambda\in\Lambda (W)} \subset\Sigma W$ be a family of good idempotents and denote
 $\Sigma_{\Lambda} W$ the subspace of $\Sigma W$ generated by
$\{E_\lambda\}_{\lambda\in \Lambda (W)}$ then
\begin{enumerate}
\item $\Sigma_{\Lambda} W$ is a subalgebra of $\Sigma W$ and 
$(E_\lambda)_{\lambda\in \Lambda (W)}$ is one of its bases;
\item $\Sigma W = \ker\Phi\oplus \Sigma_{\Lambda} W$;
\item $\Phi$ is an isomorphism from $\Sigma_{\Lambda} W$ to $\Gamma W$.
\end{enumerate}
\end{lem}
\begin{proof} As the $E_\lambda$ are orthogonal idempotents, $\Sigma_{\Lambda} W$ is 
an algebra.  To prove that $(E_\lambda)_{\lambda\in \Lambda (W)}$ is a basis of $\Sigma_{\Lambda} W$, one just has to show that
 $(E_\lambda)_{\lambda\in \Lambda (W)}$ is free. Let $a_\lambda\in\mathbb Q$ such that
 $$
 \sum_{\lambda\in\Lambda (W)} a_\lambda E_\lambda = 0 .
 $$
 As $E_\lambda \in\Sigma W$, this implies
 \begin{eqnarray*}
 \Phi(\sum_{\lambda\in\Lambda (W)} a_\lambda E_\lambda) &=& \sum_{\lambda\in\Lambda (W)} a_\lambda \Phi(E_\lambda ) \\
 &=& \sum_{\lambda\in\Lambda (W)} a_\lambda \xi_\lambda \\
 &=& 0 .
 \end{eqnarray*}
 As $(\xi_\lambda)_{\lambda\in\Lambda (W)}$ is a basis of $\Gamma W$, $a_\lambda = 0$ for all $\lambda\in\Lambda (W)$. Therefore,
 $(E_\lambda)_{\lambda\in \Lambda (W)}$ is a basis of $\Sigma_{\Lambda} W$ and $\dim\Sigma_{\Lambda} W = \card{\Lambda (W)}$.
 
Let $x\in \ker\Phi\cap \Sigma_{\Lambda} W$, 
then $\Phi(x)=\sum_{\lambda\in\Lambda (W)} a_{\lambda}\xi_\lambda =0$ where $a_\lambda \in \mathbb Q$. 
For any $\lambda\in\Lambda( W)$ and $w$ a $\lambda$-element, $\Phi(x)(w)=a_\lambda =0$. 
Therefore $x=0$. This implies  $\ker\Phi\cap \Sigma_{\lambda} W=\{0\}$ thus $\ker\Phi\oplus \Sigma_{\Lambda} W$. 
As $\dim\Sigma_{\Lambda} W + \dim \ker\Phi = \dim \Sigma W$,  $\Sigma W = \ker\Phi\oplus \Sigma_{\Lambda} W$. 

Therefore, the restriction of $\Phi$ to $\Sigma_{\lambda} W$ is an isomorphism 
from $\Sigma_{\lambda} W$ to $\Gamma W$.
\end{proof}

 \begin{cor} Any  family of good idempotents  is a family of othogonal idempotents. 
 \end{cor}
 \begin{proof} Let $\lambda,\mu$ be two distincts Coxeter classes. Then $\Phi(E_\lambda E_\mu) = \xi_\lambda \xi_\mu=0$, hence $E_\lambda E_\mu \in \ker\Phi$.
  As $E_\lambda E_\mu \in\Sigma_{\Lambda} W$, one has $E_\lambda E_\mu = 0$ by Lemma~\ref{hom1} $(2)$.
 \end{proof}
 \paragraph{Remark.} In \cite{bergeron}, the authors have obtained a set of good idempotents $(E_{\lambda})_{\lambda\in\Lambda (W)} 
 \subset\Sigma W$ such that  their sum for $\lambda\in\Lambda(W)$ equals the identity of $\ga W$. In this case,
 $\ga W = \bigoplus_{\lambda \in \Lambda (W)} E_\lambda \ga W$.

\begin{proof}[Proof of Theorem~\ref{mainidem}.] One may assume that $x=x_I$ and $y=x_J$ for
 some $I,J\subset S$, by Proposition~\ref{propsym}.

One has $\arrowvert C(\lambda)\cap X_I \arrowvert= \xi_\lambda (x_I)=\Phi(E_\lambda)(x_I)$. By Lemma~\ref{idem2}, 
  $\Phi(x_I)(E\lambda)=\ind I (E_\lambda) = <\chi_{\lambda},\ind I>_W$, hence ($1$) implies ($2$). \\
($2$) implies ($1$): As 
$\{E_\lambda\}_{\lambda\in \Lambda W}$ is a basis of $\Sigma_{\Lambda} W$ (Lemma~\ref{hom1}) and $\{x_I\}_
{I\subset S}$ is a basis of $\Sigma W$, one just has to show
 that $\ind I (E_\lambda) =\Phi(E_\lambda)(x_I)$. One has
 \begin{eqnarray*}
 \xi_\lambda (x_I) &=& \Phi(E_\lambda )(x_I)\\
 &=& \card{X_I \cap C(\lambda)}\\
 &=& <\chi_{\lambda},1_{W_I}^W>_W \\
 &=& \ind I (E_\lambda ) ,\qquad\textrm{by Lemma~\ref{idem2}.}
 \end{eqnarray*}
\end{proof}

\subsection{The double coset conjecture}

Let $I,J\subset S$ and $b\in X_{IJ}$. We consider the set
$$
W(I,J,b) = \set{w\in W \st w^J b^{-1} = (wb^{-1})_{I}} ,
$$
where $(w^{J} , w_J)$ (resp. $((wb^{-1})^{I} , (wb^{-1})_{I} )$) are the parabolic components of $w$ (resp . $wb^{-1}$). 

In \cite{hohlweg}, the authors have shown that $W(I,J,e)=W(J,I,e)$, for all $I,J\subset S$. For $b\not = e$, $W(I,J,b)$ is not generally equal to
$W(J,I,b^{-1})$. They stated the following conjecture.

\begin{conj}\label{conj2} If $I,J\subset S$ then for all $b\in X_{IJ}$, $\card{W(I,J,b)}=\card{W(J,I,b^{-1})}$.
\end{conj}

 The following proposition shows that Conjecture~\ref{conj2} implies Conjecture~\ref{conj1}.
 
\begin{prop}\label{wset} Let $I,J\subset S$ then
$$
\Phi(x_J)(x_I) = \sum_{b\in X_{IJ}} \card{W(I,J,b)}  .
$$
\end{prop}

\begin{lem}\label{lemdcc} Let $I,J\subset S$ and $b\in X_{IJ}$. Let $w\in W$ then
\begin{enumerate}
\item\label{1r} if $w \in X_{J\cap b^{-1} I b}$, $w_J \in X^{J}_{J\cap b^{-1} I b}$;
\item\label{dcc1} If $w\in W(I,J,b)$,  $w\in W_I b W_J$ and $w^{J} \in X^{I}_{I\cap bJb^{-1}}b$;
\item  $w\in W(I,J,b)$ if and only if there are $a\in X_{I\cap bJb^{-1}}^{I}$, $c\in X^{J}_{J\cap b^{-1}Ib}$ and
$u\in X_I$ such that $w = abc = uab$. Moreover, $a,c,u$ are unique.
\end{enumerate}
\end{lem}
\begin{proof}
\begin{enumerate}
\item As $X_J X^{J}_{J\cap b^{-1}Ib} =X_{J\cap b^{-1}Ib}$, $w_J \in X^{J}_{J\cap b^{-1} I b}$. 
\item One has
$$
W_I w W_J = W_I w^{J}w_J W_J = W_I w^{J} W_J = W_I (wb^{-1})_I b W_J = W_I b W_J .
$$
As seen in Section~\ref{secalg}, $w^{J} \in X^{I}_{I\cap bJb^{-1}}b$.
\item One has $w^{J} = ab$ with a unique $a\in X_{I\cap bJb^{-1}}^{I}$, by Lemma~\ref{dcc1}. 
Let $u=w^{I}\in X_I$ and $c=w_{J}$ uniquely determined as parabolic components. Then 
\begin{eqnarray*}
w=uab &\in& X_I X_{I\cap bJb^{-1}}^{I} b\\
&=& X_{I\cap bJb^{-1}}b,\qquad \textrm{see Section~\ref{secalg}}\\
&=& X_{J\cap b^{-1}Jb},\qquad \textrm{see Section~\ref{secalg}}.
\end{eqnarray*}
By (\ref{1r}),  $c=\in X^{J}_{J\cap b^{-1}Ib}$. 

If they are $a\in X_{I\cap bJb^{-1}}^{I}$, $c\in X^{J}_{J\cap b^{-1}Ib}$ and
$u\in X_I$ such that $w = abc = uab$, then $w^{J}b^{-1} = abb^{-1}=a= (wb^{-1})_I $. Thus $w\in W(I,J,b)$.
\end{enumerate}
\end{proof}

\begin{proof}[Proof of Proposition~\ref{wset}] As 
$
1_{W_J}^{W} (u) = \card{\set{g\in X_J \st ug \in gW_J }}
$ (see \cite{humphreys}, Section $3.14$), 
\begin{eqnarray*}
\Phi(x_J)(x_I) &=& \sum_{u\in X_{I}} 1_{W_J}^{W} (u)\\
&=& \card{\set{(u,g)\in X_I \times X_J \st ug \in gW_J}}  .
\end{eqnarray*}
As $X_J$ is the disjoint union of $X^{I}_{I \cap bJb^{-1}}b$ for all $b\in X_{IJ}$, the double parabolic components of $g\in X^{I}_{I \cap bJb^{-1}}b$
 are $(a , b, e)$. Hence,
\begin{eqnarray*}
\Phi(x_J)(x_I) &=& \sum_{b\in X_{IJ}} \card{\set{(u,a)\in X_I \times X^{I}_{I \cap bJb^{-1}} \st uab \in ab W_J}}\\
&=&\sum_{b\in X_{IJ}} \card{W' (I,J,b)} .
\end{eqnarray*}
One just has to show that $\card{W(I,J,b)}=\card{W' (I,J,b)}$. 

By Lemma~\ref{lemdcc}, ($3$), for any $w\in W$ one has unique $a\in X_{I\cap bJb^{-1}}^{I}$, $c\in X^{J}_{J\cap b^{-1}Ib}$ and
$u\in X_I$ such that $w = abc = uab$. Therefore, the following mapping
\begin{eqnarray*}
f\, :\ W(I,J,b) &\longrightarrow & W'(I,J,b) \\
w & \longmapsto& (u , a )
\end{eqnarray*}
is well defined and is injective  (parabolic components are unique and determine uniquely each $w\in W$). 

Let $(w,a)\in W'(I,J,b)$, then 
$$
(wab)^{J} b^{-1} = (ab) b^{-1} = a = ((wab)b^{-1})_I .
$$
Then $wab \in W(I,J,b)$ and $f\, (wab) = (w,a)$. Therefore $f$ is a bijection and the proposition is proved. 
\end{proof}

\subsection{Results}

The following proposition proves Conjecture~\ref{conj2} for any $I\subset S$ and $J = \emptyset , S$.

\begin{prop}\label{dcc2} If $(W,S)$ is a finite Coxeter system and if $I\subset S$, then
\begin{enumerate}
\item $\card{W(I,\emptyset , b)} =\card{W(\emptyset ,I , b^{-1})} = \card{W_I}$, for all $b\in X_{I\emptyset}$;
\item $X_{IS}=\set{e}$ and $W(I,S , e) =W(S ,I , e)= X_I$.
\end{enumerate}
\end{prop}

\begin{proof}
\begin{enumerate}
\item Observe that $X_{\emptyset I} = X_I$. If $w\in W$ then $(w^{\emptyset},w_{\emptyset} )= (w,e)$. Let $b\in X_{I}^{-1}$. Then, on one hand,
\begin{eqnarray*}
W(I,\emptyset , b)&=&\set{w\in W\st w^{\emptyset}b^{-1} = (wb^{-1})_I}\\
&=&\set{w\in W\st w = (wb^{-1})_I b}\\
&=&\set{w\in W\st w \in W_Ib}\\
&=& W_I b
\end{eqnarray*}
On the other hand,
\begin{eqnarray*}
W(\emptyset ,I, b^{-1})&=&\set{w\in W\st w^{I}b = (wb)_\emptyset}\\
&=&\set{w\in W\st w^I b= e }\\
&=&\set{w\in W\st w^I = b^{-1} }\\
&=&\set{w\in W\st w \in b^{-1}W_I}\\
&=& b^{-1} W_I 
\end{eqnarray*}
\item Observe that $X_{IS} = \set{e}$ then $(w^{S},w_{S} )= (e,w)$. On one hand,
\begin{eqnarray*}
W(I,S, b)&=&W(I,S,e)\\
&=&\set{w\in W\st w^{S} = w_I}\\
&=&\set{w\in W\st e= w_I }\\
&=&\set{w\in W\st w^I \in X_I }\\
&=& X_I 
\end{eqnarray*}
On the other hand,
\begin{eqnarray*}
W(S,I, b^{-1})&=&W(S,I,e)\\
&=&\set{w\in W\st w^{I} = w_S}\\
&=&\set{w\in W\st W^{I}= w}\\
&=& X_I 
\end{eqnarray*}
\end{enumerate}
\end{proof}

The following corollary, that proves Conjecture~\ref{conj1} for any $I\subset S$ 
and $J = \emptyset , S$, may be proved directly using Proposition~\ref{propsym}.

\begin{cor}\label{cordcc1} Let $(W,S)$ be a finite Coxeter system then for any $I\subset S$
\begin{enumerate}
\item $\Phi(x_\emptyset)(x_I) =\Phi(x_I)(x_\emptyset)$,
\item $\Phi(x_S)(x_I) =\Phi(x_I)(x_S)$.
\end{enumerate}
\end{cor}
\begin{proof} It is an immediate consequence of Proposition~\ref{dcc2} and Propostion~\ref{wset}.
\end{proof}

\begin{prop}\label{dcc3} Let $(W,S)$ be a finite Coxeter system, $I,J\subset S$ and $b\in X_{IJ}$ such that $bJb^{-1} \subset I$. Then
$$
\card{W(I,J,b)} = \card{W(J,I,b^{-1})} = \frac{\card{W_I}}{\card{W_J}} .
$$
\end{prop}
\begin{proof}  On one hand, observe that $X^{J}_{J\cap b^{-1}Ib} = X^{J}_J 
=\set{e}$. Thus $W(I,J,b) = X^{I}_{I\cap bJb^{-1}}b$, by Lemma~\ref{lemdcc}. 
Therefore,
\begin{eqnarray*}
\card{W(I,J,b)} &=& \card{X^{I}_{I\cap bJb^{-1}}}\\&=&\frac{\card{W_I}}{\card{W_{I\cap b J b^{-1}}}}\\& =&
\frac{\card{W_I}}{\card{W_I \cap b W_J b^{-1}}}\\& =& \frac{\card{W_I}}{\card{ b W_J b^{-1}}} \\&=& \frac{\card{W_I}}{\card{  W_J }}\, , \qquad\textrm{as }
 bW_J b^{-1} \subset W_I\ .
\end{eqnarray*}

On the other hand, $W(J,I,b^{-1}) \subset b^{-1}X^{I}_{I\cap bJb^{-1}}$ by Lemma~\ref{lemdcc}. Let $a\in X^{I}_{I\cap bJb^{-1}}$, then
\begin{eqnarray*}
b^{-1}ab &\in & X_I X^{I}_{I\cap bJb^{-1}} b\\
&=& X_{J\cap b^{-1}Ib}\\
&=& X_J,\ \textrm{as } J\subset b^{-1} I b\ .
\end{eqnarray*}
Hence $b^{-1} a = ua b \in W(J,I,b^{-1})$ by Lemma~\ref{lemdcc}. Therefore $W(J,I,b^{-1}) = b^{-1}X^{I}_{I\cap bJb^{-1}}$ and the lemma is proved.
\end{proof}

\subsection{The case of a single generator}

We denote $t$ for $\set{t}\subset S$ to simplify the notations.

\begin{thm}\label{dccsing} Let $(W,S)$ be a finite Coxeter system, $t\in S$ and $I\subset S$ then for all $b\in X_{It}$
\begin{eqnarray*}
\card{W(I,t , b)} 
           &=& \card{W(t , I, b^{-1})}   .
\end{eqnarray*}
\end{thm}

\begin{cor}\label{died} The double coset conjecture holds for dihedral groups.
\end{cor}
\begin{proof}
 As for $I_2 (m)$, $S=\set{s,t}$, the subsets of $S$ are $\emptyset$, ${s}$, ${t}$, $S$, one concludes 
with Corollary~\ref{cordcc1} and Theorem~\ref{dccsing}.
\end{proof}

\begin{cor}\label{cordcc2} Let $(W,S)$ be a finite Coxeter system. Then 
$$
\Phi(x_I)(x_t)= \Phi(x_t)(x_I)  ,
$$
for any $t\in S$ and any $I\subset S$. On other words,
$$
\sum_{w \in X_I} \ind t (w) = \sum_{g \in X_t} \ind I (g), 
$$
for any $t\in S$ and any $I\subset S$.
\end{cor}

\begin{proof}[Proof of the theorem] 
\begin{enumerate}
\item If $btb^{-1} \in I$, one concludes by Proposition~\ref{dcc3}.
\item From now on, one supposes that $btb^{-1} \notin I$. Thus
$$
X_{I\cap btb^{-1}}^{I} = W_I \Longrightarrow  W_I b\subset X_t
$$
and
$$
X_{t \cap b^{-1}Ib}^{t} = W_{t} \Longrightarrow W_t b^{-1}\subset X_I
 .
$$
\item Let $d\in X_{II}$ such that $btb \in W_I d W_I$. Then $d^{2}=e$. Indeed, as $d\in X_{II}$
 is the unique element of minimal length in $W_I d W_I$ and as 
$$
(btb^{-1})^{-1} = btb^{-1}\in W_I d^{-1} W_I = W_I d W_I \ ,
$$
 one has $d=d^{-1}$.
\item There is a unique $a_0 \in X_{I\cap dId}^{I}$ such that $btb^{-1} = a_0 d a_{0}^{-1}$.
Indeed, there is a unique $(a_0 ,a_1 )\in X^{I}_{I\cap dId}\times W_I$ such that
$btb^{-1} = a_0 d a_1$. Thus 
$$
tb^{-1} a_1 = b^{-1} a_0 d \in X_I X^{I}_{I\cap dId}d = X_{I\cap dId}d = X_{I\cap dId} .
$$
Then $a_1 \in X_{I\cap dId}^{I}$, by Lemma~\ref{lemdcc} $(1)$.

Finally one has 
$$
btb^{-1} = (btb^{-1} )^{-1}=a_o d a^{-1}_1 = a_1 d a^{-1}_{0} . 
$$
By unicity of  double parabolic components, $a_0 = a_1$.
\item\label{formule1}  One defines a mapping:
\begin{eqnarray*}
\sigma \, :\ W_I &\longrightarrow & btb^{-1} W_I btb^{-1} \\
a & \longmapsto& btb^{-1}abtb^{-1}
\end{eqnarray*}
As $\sigma$ is the restriction of the conjugation by $btb^{-1}$ on $W_I$, it is a isomorphism.
Denote
$$
Fix_{\sigma} W_I =\set{a\in W_I \st \sigma (a) = a} .
$$
As $e\in W_I$ and $\sigma (e) = e$, this set is non empty. 
One will  show that
$$
\card{W(I,t , b)}=\card{W_I}+ \card{ Fix_{\sigma} W_I} .
$$
In fact, one will prove that
$$
(\star)\qquad W(I,t , b)=W_I b \amalg a^{-1}_0 Fix_{\sigma} W_I bt .
$$
It is immediate that $W_I b \subset W(I,t,b)$.\\
 Let $w = a_{0}^{-1}a bt \in Fix_\sigma W_I$.
 As $a_{0}^{-1}a \in W_I$, $a_{0}^{-1}a b \in X_t$, by $2$. Moreover,
\begin{eqnarray*}
 a_{0}^{-1}a b t b^{-1}a^{-1}a_0 &=& a_{0}^{-1}(a a_0 d a^{-1}_0 )a^{-1}a_0 ,\ \textrm{as }
 btb^{-1}=a_0 d a^{-1}_0 \\
&=& a_{0}^{-1}(a_0 d a^{-1}_0  a) a^{-1}a_0 ,\ \textrm{as }
 a\in Fix_\sigma W_I \\
&=&d\in X_{II}\subset X_I  .
\end{eqnarray*}
 Therefore $w\in W(I,t,b)$ that implies $a^{-1}_0 Fix_{\sigma} W_I bt\subset W(I,t,b)$.

Conversely, let $w\in W(I,t,b)$, then there is $a\in W_I$ such that $w=ab\in W_I b$ or $w=abt$. 
In the second case $abtb^{-1}a^{-1}=u \in X_I$, by Lemma~\ref{lemdcc}. As $u$ is an involution,
 $u\in X_{II}$. But $u \in aW_I d W_I a^{-1} = W_I d W_I$ therefore $u=d$.\\
One writes $w=a_{0}^{-1} (a_0 a) bt$. Denote $a' = a_0 a \in W_I$ thus
$$
a_{0}^{-1} a' btb^{-1} a'^{-1}a_0=abtb^{-1}a^{-1}=d .
$$
therefore
$$
a' btb^{-1} a'^{-1}=a_0 d a^{-1}_0 = btb^{-1} .
$$
This implies $a'\in Fix_\sigma W_I$ thus $w\in a^{-1}_0 Fix_{\sigma} W_I bt$ and 
$(\star)$ is proved.

\item\label{6} If $a\in W_I $ and  $d_a \in X_{t\set{t}}$ such that $b^{-1} a b \in W_{t} d_a W_{t}$.
 Then the following propositions are equivalent: 
 \begin{enumerate}
\item[i)] $td_a = d_a t$,
\item[ii)] $b^{-1} a b = tb^{-1} a b t = d_a $,
\item[iii)] $\set{b^{-1} a ,tb^{-1} a}\subset W(t , I, b^{-1})$,
\item[iv)] $b^{-1} a b = tb^{-1} a b t$.
\end{enumerate}
$i)\Rightarrow ii)$: If $t d_a = d_a t$ then $X_{t\cap d_a t d_{a}^{-1}}^t =\set{e}$ and $W_t d_a W_t =\set{d_a , d_a t}$. As $b^{-1} a b \in W_t d_a W_t$, 
$b^{-1} a b = d_a$ or $b^{-1} a b = d_a t$. In the second case, that is $ab= bd_a t$, one has $bd_a \in X_t d_a =X_{d^{-1}_a t d_a} = X_t$. Thus $ab=bd_a t \notin X_t$, which is a contradiction with $ab\in W_I b \subset X_t$.
 Therefore $b^{-1} a b = d_a = td_a t = tb^{-1} a b t$. \\
$ii)\Rightarrow iii)$: If $b^{-1} a b = tb^{-1} a b t = d_a$, then $b^{-1} a = d_a b$ and $tb^{-1} a = d_a t b$. Therefore, $\set{b^{-1} a ,tb^{-1} a}
\subset W(t , I, b^{-1})$, by Lemma~\ref{lemdcc} $3$.\\
$iii)\Rightarrow i)$: If $\set{b^{-1} a ,tb^{-1} a}\subset W(t , I, b^{-1})$, thus $b^{-1} ab$ and $tb^{-1} abt$ are in  $X_t$.
One supposes that $t d_a \not = d_a t$. Then $X^{t}_{t\cap d_a t d^{-1}_a} = W_t$, hence $W_t d_a \subset X_t$.\\
Therefore, 
$$
b^{-1} ab =\left\{\begin{array}{cc} d_a \in X_t\ &\textbf{ or}\\ td_a \in X_t &\end{array}\right. .
$$
Thus
$$
tb^{-1} abt =\left\{\begin{array}{cc} td_a t\notin X_t\ &\textbf{ or}\\ d_a t\notin X_t &\end{array}\right. .
$$
Which is a contradiction with $iii)$. Hence $b^{-1} ab = tb^{-1} abt \in X_t$. 
Therefore, 
$$
b^{-1} ab = tb^{-1} abt =\left\{\begin{array}{cc} d_a \in X_t\ &\textbf{ or}\\ td_a \in X_t &\end{array}\right. .
$$
Finally, either $b^{-1} ab = d_a$, then $tb^{-1} ab t= td_a t = d_a$ hence $d_a t = t d_a $; or
 $b^{-1} ab = td_a$, then $tb^{-1} ab t= d_a t = td_a$ hence $d_a t = t d_a $.\\
$i)\iff iv)$: As $i)\Rightarrow ii)$, one has in particular $b^{-1} a b = tb^{-1} a b t$ thus $i)\Rightarrow iv)$. 

If $b^{-1} a b = tb^{-1} a b t$, so
$$
b^{-1} ab = tb^{-1} abt = \left\{\begin{array}{cccc} d_a = td_a t\in X_I \ &\textbf{ or}\\ td_a = d_a t \in X_I \ &\textbf{ or}\\ d_a t = td_a \notin X_I 
\ &\textbf{ or}\\td_a t = d_a \notin X_I &\end{array}\right. .
$$
Observe that $td_a = d_a t$ in all the cases.

\item\label{formule2} One has:
$$
\card{W(I,t,b^{-1})} = \card{W_I} +
 \card{\set{a \in W_I \st d_a t = t d_a }} .
 $$
  Indeed, let $a\in W_I$ such that $d_a t \not = t 
d_a$ thus $X^{t}_{t\cap d_a t d^{-1}_a} = W_t$. As $b^{-1} ab \in W_t d_a W_t$, one has
$$
b^{-1} ab = \left\{\begin{array}{cc} d_a \in X_I \ &\textbf{ or}\\ td_a \in X_I \ &\textbf{ or}\\ d_a t \notin X_I 
\ &\textbf{ or}\\td_a t \notin X_I &\end{array}\right. .
$$
Observe that either
$$
b^{-1} a \in W(t,I,b^{-1})\ \textrm{ and }\ tb^{-1} a \notin W(t,I,b^{-1})
$$
or
$$
b^{-1} a \notin W(t,I,b^{-1})\ \textrm{ and }\ tb^{-1} a \in W(t,I,b^{-1}) .
$$
As $W_I$ is the disjoint union of $\set{a\in W_I \st d_a t = t d_a}$ and $\set{a\in W_I \st d_a t 
\not = t d_a}$, one has:
\begin{eqnarray*}
\card{W(t,I,b^{-1})}&=& \sum_{a\in W_I \st d_a t = t d_a} 2 + \sum_{a\in W_I \st d_a t 
\not = t d_a} 1\, , \ \textrm{ by (\ref{6})}\\
&=& 2 \card{\set{a \in W_I \st d_a t = t d_a }} + \card{\set{a \in W_I \st d_a t \not= t d_a }}\\
&=& \card{W_I} +  \card{\set{a \in W_I \st d_a t = t d_a } }  .
\end{eqnarray*}

\item By \ref{formule1} and \ref{formule2}, one just has to show that $\card{Fix_\sigma W_I}
 = \card{\set{w \in W_I \st d_a t = t d_a }}$ to finish the proof. 
\begin{eqnarray*}
a\in Fix_\sigma W_I &\iff& abtb^{-1} = btb^{-1} a\\
&\iff& b^{-1}ab = tb^{-1}ab t\\
&\iff& d_a t = t d_a\ , \qquad\textrm{by \ref{6}}\\
&\iff& a\in\set{w \in W_I \st d_a t = t d_a }
\end{eqnarray*}
\end{enumerate}
\end{proof}

\subsection{Proof of Theorem~\ref{conjtrue}}

By Proposition~\ref{propsym} and Proposition~\ref{irred}, one shows the theorem by using the classification of
 irreducible finite Coxeter groups (see \cite{humphreys}).

\paragraph{Case $A_n$:} This case has been shown by J\" ollenbeck and Reutenauer in \cite{jollen}. 
An idea of their proof is to be found in the remark following Theorem~\ref{mainidem}.

\paragraph{Case $I_2 (m)$:} It is given by Proposition~\ref{wset} and Corollary~\ref{died}.

\paragraph{Cases $E_6$, $E_7$, $E_8$, $F_4$, $H_3$, $H_4$:} We begin to give a lemma.

\begin{lem}\label{1isim1j} Let $I,J\subset S$, then
$$
\sum_{w\in X_I} 1^{W}_{W_J} (w) = \sum_{C\in Cl(W)} \card{X_I \cap C} 1_{W_J}^{W} (w_C)\ ,
$$
where $w_C$ is an arbitrary element in $C$.
\end{lem}
\begin{proof} One has
$$
\sum_{w\in X_I} 1^{W}_{W_J} (w)=\sum_{C\in Cl(W)} \sum_{w\in X_I \cap C} 1^{W}_{W_J} (w) .
$$
As $1^{W}_{W_J}$ is constant in conjugacy classes, one chooses an arbitrary element  $w_C \in C$. Hence
\begin{eqnarray*}
\sum_{w\in X_I} 1^{W}_{W_J} (w)&=&\sum_{C\in Cl(W)} \sum_{w\in X_I \cap C} 1^{W}_{W_J} (w_C)\\
&=&\sum_{C\in Cl(W)}1^{W}_{W_J} (w_C)\sum_{w\in X_I \cap C} 1\\
&=&\sum_{C\in Cl(W)} \card{X_I \cap C} 1_{W_J}^{W} (w_C)  .
\end{eqnarray*}
\end{proof}

As, if $J\sim J'$, $1^{W}_{W_J}=1^{W}_{W_J'}$  and $\card{X_I \cap C} =\card{X_J \cap C}$ 
for some conjugacy class $C$ (see Proposition~\ref{propker}), it is sufficient to consider 
Coxeter Classes of $I$ and $J$ in Lemma~\ref{1isim1j}. Therefore, we construct two matrices 
$$
A = (a_{\lambda C})_{\lambda \in \Lambda(W) , C \in Cl(W)}
\textrm{ and   }
 B=(b_{C' \mu })_{C' \in Cl(W) , \mu\in \Lambda(W)},
 $$
  where
$
a_{\lambda C} = \card{X_I \cap C}, \ \textrm{for some } I\in \lambda 
$
and
$
 b_{C' \mu }= \varphi_\mu (w_C') 
 $.
Then $AB =D= (d_{\lambda \mu})_{\lambda , \mu \in \Lambda(W)}$ where
$$
d_{\lambda \mu} = \sum_{C\in Cl(W)} \card{X_I \cap C} 1_{W_J}^{W} (w_C)= \sum_{w\in X_I} 1^{W}_{W_J} (w),\ \textrm{for some } I\in \lambda,\ J \in \mu .
$$
Hence to prove Conjecture~\ref{conj1} amounts to prove that $D$ is a symmetric matrix.\\

For these exceptional cases,
 we used the GAP part of CHEVIE (\cite{chevie} and \cite{gap}) 
to obtain $D$. This program contains some useful functions to compute the matrix $A$ and the matrix $B$, therefore the matrix $D$.
\begin{enumerate}
\item the GAP functions {\tt CharTable}, {\tt PositionId}, {\tt StoreFusion} and {\tt Induced} are used to compute the coefficients of $A$;
\item The functions {\tt ReflectionSubgroup}, {\tt ReducedRightCosetRepresentatives} and {\tt PositionClass} are used to compute the coefficients of $B$.
\end{enumerate}

We just have to give a  set of representants of Coxeter classes. Following the notation in \cite{geck}, 
 the following algorithm  determines a set of representants for 
each Coxeter Class. 

\begin{algo}\label{acc} [Construction of Representant of Coxeter Classes] Given a finite Coxeter System $(W,S)$  
$\Lambda(W)_{rep}$ - a set of representant of $\Lambda(W)$ - is determinated. As  $I\sim J$ iff $\card I = \card J$, one just has to consider
 the subset of $S$ of cardinality $i$, for $i=1,\dots , \card{S}$ to construct $\Lambda(W)_{rep}$.
{\rm
\begin{enumerate}
\item \textbf{[Initialize]} Let $PS_i = \set{I_1 , \dots ,I_m }$ be the set of subsets of $S$ of cardinality $i$. 
Set $\Lambda(W)^{i}_{rep} \leftarrow I_1$ and 
$$
E_1 \leftarrow \set{I_p \st \ c_{I_p}\sim c_{I_1} } .
$$
\item \textbf{[Next $k$]} If $I_k \in E_{k-1}$, then set $k\leftarrow k+1$, $E_{k+1} \leftarrow E_k$. Otherwise set 
$$
\Lambda(W)^{i}_{rep} \leftarrow \Lambda(W)^{i}_{rep} \cup I_k
$$
and
$$
E_k \leftarrow E_{k-1} \cup \set{I_p \st \ c_{I_p}\sim c_{I_k} }
$$
and $k\leftarrow k+1$.\\
At this step, each subset of $S$ in $\Lambda(W)^{i}_{rep}$ is not conjugated to each of the others and  each subset of $S$ in $E_k$ has a representant 
of it Coxeter class in $\Lambda(W)^{i}_{rep}$.
\item \textbf{[Loop]} If $k < m+1$ repeat step $2$ else stop.
\end{enumerate}
}
Then $$\Lambda(W)_{rep} = \bigcup_{i=0}^{\card{S}} \Lambda(W)^{i}_{rep} .$$
\end{algo}

By 
Corollary~\ref{cordcc2}, one just has to take the submatrix $A'$ of $A$ and $B'$ of $B$
 indexed by Coxeter Classes $\lambda , \mu$ such that  
$\card I , \card J \geq 2$ for some $I\in \lambda$ and $J\in \mu$. The above algorithm can be initialize to $i=2$. 

As  $\card{E_8}> 6.10^{8}$, using Corollary~\ref{cordcc2} is necessary to restrain 
the cardinality of cosets in this very large case. Finally, still in this case, we cannot directly construct and stock the set of minimal 
coset representatives $X_I$ with GAP. Indeed for a set of two commuting generators in $E_8$,
 $\card{X_I} > 150.10^{6}$. So we used the fact that for $I\subset S$ and $s\in S-I$ then
 $$
 X_I = X_K X^{K}_I\ ,
 $$
where $K = S-\set s$.  If we keep in memory $X_I$, it take the memory of 
$X_K$ times $X^{K}_I$, rather than if we keep in memory $X_K$ and $X^{K}_I$, 
they take the memory for $X_K$ plus $X^{K}_I$.\\

Finally, we just have to compute the submatrix $D'$ of $D$, indexed by Coxeter Classes $\lambda , \mu$ such that  
$\card I , \card J \geq 2$. Observe in the following that $D'$ is a symmetric matrix.
(The following matrix and tables have been obtained by one week of computation with GAP-CHEVIE; the original program
 also verifyied that $D'$ is symmetric).  This complete our proof.

\paragraph{Remark about GAP-CHEVIE notations:} The set $S$ of generators are given in CHEVIE by integers from $1$ to $n$, 
see \cite{chevie} for more informations. The set of representants \{I\} of Coxeter classes is given as a set of Coxeter $I$--elements.

\paragraph{Coxeter group of type $H_3$}
The set of representants of Coxeter classes  that we consider is:
$$
\left\{ \begin{array}{cccc} 
12,&123,&13,&23\end{array} \right\} .
$$
Therefore the submatrix of $D'$ indexed by them is
$$
\left( \begin{array}{llll} 
24&12&40&32\\
12&1&30&20\\
40&30&52&46\\
32&20&46&38\\
\end{array} \right) 
$$

\paragraph{Coxeter group of type $H_4$}
The set of representants of Coxeter classes that we consider is
$$
\left\{ \begin{array}{cccccccc} 
12,&123,&1234,&124,&13,&134,&23,&234\end{array} \right\} .
$$
Therefore the submatrix of $D'$ indexed by them is
$$
\left( \begin{array}{cccccccc} 
4080&2040&1440&2840&5760&3280&4800&2640\\
2040&404&120&1116&4112&1602&2952&958\\
1440&120&1&720&3600&1200&2400&600\\
2840&1116&720&1828&4744&2288&3672&1664\\
5760&4112&3600&4744&7144&5104&6360&4592\\
3280&1602&1200&2288&5104&2730&4080&2132\\
4800&2952&2400&3672&6360&4080&5484&3504\\
2640&958&600&1664&4592&2132&3504&1510\\
\end{array} \right) 
$$

\paragraph{Coxeter group of type $F_4$}
The set of representants of Coxeter classes that we consider is:
$$
\left\{ \begin{array}{ccccccccc} 
12,&123,&1234,&124,&13,&134,&23,&234,&34\end{array} \right\} .
$$
Therefore the submatrix of $D'$ indexed by them is:
$$
\left( \begin{array}{ccccccccc} 
396&252&192&324&552&396&456&312&600\\
252&70&24&148&364&170&232&84&312\\
192&24&1&96&288&96&144&24&192\\
324&148&96&224&432&250&312&170&396\\
552&364&288&432&600&432&512&364&552\\
396&170&96&250&432&224&312&148&324\\
456&232&144&312&512&312&416&232&456\\
312&84&24&170&364&148&232&70&252\\
600&312&192&396&552&324&456&252&396\\
\end{array} \right) .
$$

\paragraph{Coxeter group of type $E_6$}
The set of representants of Coxeter classes that we consider is
$$
\left\{ 
12,\, 123,\, 1234,\, 12345,\, 123456,\, 12346,\, 1235,\, 12356,\, 1245,\, 125,\, 13,\, 134,\, 13456,\, 1356,\, 2345 \right\} .
$$
Therefore the submatrix of $D'$ indexed by them is:
\begin{center}
\tiny{$$
\left( \begin{array}{ccccccccccccccc} 
29136&21744&16128&13728&12960&14544&18144&15456&17088&23040&26784&19776&14448&17856&15408\\
21744&12858&6882&4890&4320&5616&9252&6564&7950&14712&18432&10488&5412&8640&6300\\
16128&6882&1704&642&432&1088&3976&1864&2648&9108&12024&4432&840&3144&1380\\
13728&4890&642&87&27&341&2540&926&1378&7092&9432&2632&162&1776&444\\
12960&4320&432&27&1&216&2160&720&1080&6480&8640&2160&72&1440&270\\
14544&5616&1088&341&216&661&3080&1304&1882&7800&10368&3336&464&2304&852\\
18144&9252&3976&2540&2160&3080&6152&3894&4908&11256&14472&6976&2878&5460&3576\\
15456&6564&1864&926&720&1304&3894&1996&2690&8700&11448&4320&1104&3132&1602\\
17088&7950&2648&1378&1080&1882&4908&2690&3604&10080&13176&5552&1652&4140&2274\\
23040&14712&9108&7092&6480&7800&11256&8700&10080&16392&20016&12576&7644&10752&8544\\
26784&18432&12024&9432&8640&10368&14472&11448&13176&20016&23832&15984&10224&14112&11160\\
19776&10488&4432&2632&2160&3336&6976&4320&5552&12576&15984&7888&3088&6264&3864\\
14448&5412&840&162&72&464&2878&1104&1652&7644&10224&3088&250&2052&624\\
17856&8640&3144&1776&1440&2304&5460&3132&4140&10752&14112&6264&2052&4644&2784\\
15408&6300&1380&444&270&852&3576&1602&2274&8544&11160&3864&624&2784&1074\\
\end{array} \right) .
$$}
\end{center}

\paragraph{Coxeter group of type $E_7$}
The set of representants of Coxeter classes that we consider is\\
$
\{ 
\lambda_1=12,\ \lambda_2=123,\ \lambda_3=1234,\ \lambda_4=12345,\ \lambda_5=123456,\ \lambda_6=1234567,\\ \lambda_7=123457,\ \lambda_8=12346, \ \lambda_9=123467,\ \lambda_{10}=1235,\ 
\lambda_{11}=12356,\\ \lambda_{12}=123567,\ \lambda_{13}=12357,\ \lambda_{14}=1245,\ \lambda_{15}=124567, \ 
\lambda_{16}=12457,\\ \lambda_{17}=125,\ \lambda_{18}=1257,\ \lambda_{19}=13,\ \lambda_{20}=134,\ 
\lambda_{21}=13456,\ \lambda_{22}=134567,\\ \lambda_{23}=13467, \ \lambda_{24}=1356,\ \lambda_{25}=2345,\ \lambda_{26}=234567,\ \lambda_{27}=23457,\\ \lambda_{28}=24567,\ \lambda_{29}=2457,\ \lambda_{30}=257
 \}.
$\\
We give here  the submatrix of $D'$ indexed by this Coxeter classes. 
We have decomposed this matrix into three tables.
\begin{center}
\tiny{$$
\begin{array}{|c|c|c|c|c|c|c|c|c|c|c|}\hline 
&\lambda_{1} &\lambda_{2} &\lambda_{3} &\lambda_{4} &\lambda_{5} &\lambda_{6} &\lambda_{7} &\lambda_{8} &\lambda_{9} &\lambda_{10} \\
\hline 
\lambda_{1}&1733376&1320192&981504&823296&749568&725760&774912&877824&800256&1111296  \\       
\lambda_{2}&1320192&818544&452496&314592&258816&241920&278976&365568&300480  &601968 \\    
\lambda_{3}&981504&452496&129972&50868&29412&24192&38160&84180&50220&271872    \\ 
\lambda_{4}&823296&314592&50868&8896&2508&1512&5432&28056&11804&169328   \\  
\lambda_{5}&749568&258816&29412&2508&222&56&1330&15060&5294&131864    \\ 
\lambda_{6}&725760&241920&24192&1512&56&1&756&12096&4032&120960   \\ 
\lambda_{7}&774912&278976&38160&5432&1330&756&3198&20466&8122&145800     \\
\lambda_{8}&877824&365568&84180&28056&15060&12096&20466&52152&28862&208448   \\  
\lambda_{9}&800256&300480&50220&11804&5294&4032&8122&28862&13762&161392    \\ 
\lambda_{10}&1111296&601968&271872&169328&131864&120960&145800&208448&161392&410472   \\
\lambda_{11}&936960&427152&134928&66140&45756&40320&53716&93912&63908&261352   \\
\lambda_{12}&827904&326304&67260&22032&12334&10080&16286&41510&23056&181688  \\
\lambda_{13}&969984&463968&167688&92004&67416&60480&76680&121980&88332&295176  \\
\lambda_{14}&1044480&522096&191568&98680&68648&60480&80272&135896&94472&333824  \\
\lambda_{15}&798720&297600&46452&8768&3026&2016&5490&25720&11126&158296  \\
\lambda_{16}&933888&419808&123504&54780&35364&30240&42828&83056&53392&252776  \\
\lambda_{17}&1398528&924672&577632&440352&381312&362880&402240&489408&423648&716544  \\
\lambda_{18}&1177344&678720&347808&237840&194688&181440&210000&278736&227424&485568  \\
\lambda_{19}&1612800&1147680&761760&586080&508320&483840&535680&649440&564480&914400  \\
\lambda_{20}&1217280&689184&312000&181920&134496&120960&152448&232896&173088&472032  \\
\lambda_{21}&847872&335136&61764&14212&5636&4032&9396&35372&16456&183680  \\
\lambda_{22}&773376&276576&36156&4580&1012&576&2674&19092&7298&143496  \\
\lambda_{23}&903936&389328&100488&39404&23928&20160&30192&64996&39184&227136 \\ 
\lambda_{24}&1072128&551424&218304&122088&89688&80640&101976&159696&116328&359904  \\
\lambda_{25}&956160&427680&111312&37824&19344&15120&27168&69696&38688&252720  \\
\lambda_{26}&772608&274896&34200&3424&402&126&1806&17714&6392&141928  \\
\lambda_{27}&876288&360672&76800&21824&10096&7560&14992&45960&23512&203488  \\
\lambda_{28}&870912&351792&67284&15480&5892&4032&9900&38460&17760&194280   \\
\lambda_{29}&1067520&539280&198096&100752&69264&60480&80976&139728&96288&345264   \\
\lambda_{30}&1421568&942336&585216&443520&382464&362880&403200&494208&426240&729216  \\
\hline \end{array}
$$}
\end{center}
\begin{center}
\tiny{$$
\begin{array}{|c|c|c|c|c|c|c|c|c|c|c|}\hline 
&\lambda_{11}&\lambda_{12} &\lambda_{13} &\lambda_{14} &\lambda_{15} &\lambda_{16} &\lambda_{17} &\lambda_{18} &\lambda_{19} &\lambda_{20} \\
\hline 
\lambda_{1}  &936960&827904&969984&1044480&798720&933888&1398528&1177344&1612800&1217280 \\ 
\lambda_{2}  &427152&326304&463968&522096&297600&419808&924672&678720&1147680&689184  \\
\lambda_{3}  &134928&67260&167688&191568&46452&123504&577632&347808&761760&312000  \\
\lambda_{4}  &66140&22032&92004&98680&8768&54780&440352&237840&586080&181920  \\
\lambda_{5}  &45756&12334&67416&68648&3026&35364&381312&194688&508320&134496  \\
\lambda_{6}  &40320 &10080&60480&60480&2016&30240&362880&181440&483840&120960 \\
\lambda_{7}  &53716&16286&76680&80272&5490&42828&402240&210000&535680&152448 \\ 
\lambda_{8}  &93912&41510&121980&135896&25720&83056&489408&278736&649440&232896  \\
\lambda_{9}  &63908&23056&88332&94472&11126&53392&423648&227424&564480&173088  \\
\lambda_{10}  &261352&181688&295176&333824&158296&252776&716544&485568&914400&472032   \\
\lambda_{11}  &138040&78624&167784&188688&60868&128288&547776&332976&718560&297216 \\  
\lambda_{12}  &78624&33508&104076&113584&20118&67948&448992&248496&596160&198912  \\ 
\lambda_{13}  &167784&104076&197748&222384&84948&157728&583104&366336&756000&335808  \\
\lambda_{14}  &188688&113584&222384&254328&91056&178592&642624&410064&835200&385344  \\
\lambda_{15}  &60868&20118&84948&91056&8150&49644&421152&223440&563040&169824  \\
\lambda_{16}  &128288&67948&157728&178592&49644&117004&542112&323568&714240&288096  \\
\lambda_{17}  &547776&448992&583104&642624&421152&542112&1023360&790080&1238400&806112  \\
\lambda_{18}  &332976&248496&366336&410064&223440&323568&790080&558336&990720&554400  \\
\lambda_{19}  &718560&596160&756000&835200&563040&714240&1238400&990720&1470240&1023840  \\
\lambda_{20}  &297216&198912&335808&385344&169824&288096&806112&554400&1023840&546624  \\
\lambda_{21}  &75116&27576&102048&111768&13228&63656&460512&252768&614880&201792  \\
\lambda_{22}  &51936&15250&74904&78088&4752&41176&399936&207744&534240&149952  \\
\lambda_{23}    &108044&52608&137052&153784&35912&97240&512640&298032&679680&256704  \\
\lambda_{24}  &212328&135576&246192&281328&112752&202368&669792&434880&869760&417408  \\
\lambda_{25}  &120048&55344&152688&171888&35088&108672&555360&330480&728640&285120  \\
\lambda_{26}  &50604&14118&73320&76240&3694&39412&398688&205776&532800&147648  \\
\lambda_{27}  &88616&36032&116664&129440&19904&76672&486528&273984&645120&225792  \\
\lambda_{28}  &80400&29436&108072&120144&13788&67728&478176&263520&640800&216576  \\
\lambda_{29}  &194832&115728&228672&263712&91584&182640&660864&420240&862560&402048  \\
\lambda_{30}  &555264&451584&589824&653184&421632&546048&1042560&799488&1267200&824832  \\
\hline \end{array}
$$}
\end{center}
\begin{center}
\tiny{$$
\begin{array}{|c|c|c|c|c|c|c|c|c|c|c|}\hline 
&\lambda_{21} &\lambda_{22} &\lambda_{23} &\lambda_{24} &\lambda_{25} &\lambda_{26} &\lambda_{27} &\lambda_{28} &\lambda_{29} &\lambda_{30} \\
\hline 
\lambda_{1}  &847872&773376&903936&1072128&956160&772608&876288&870912&1067520&1421568\\
\lambda_{2}  &335136&276576&389328&551424&427680&274896&360672&351792&539280&942336\\
\lambda_{3}  &61764&36156&100488&218304&111312&34200&76800&67284&198096&585216\\
\lambda_{4}  &14212&4580&39404&122088&37824&3424&21824&15480&100752&443520\\
\lambda_{5}  &5636&1012&23928&89688&19344&402&10096&5892&69264&382464\\
\lambda_{6}  &4032 &576&20160&80640&15120&126&7560&4032&60480&362880 \\
\lambda_{7}  &9396&2674&30192&101976&27168&1806&14992&9900&80976&403200\\  
\lambda_{8}  &35372&19092&64996&159696&69696&17714&45960&38460&139728&494208\\
\lambda_{9}  &16456&7298&39184&116328&38688&6392&23512&17760&96288&426240\\
\lambda_{10}  &183680&143496&227136&359904&252720&141928&203488&194280&345264&729216\\
\lambda_{11} &75116&51936&108044&212328&120048&50604&88616&80400&194832&555264\\
 \lambda_{12}&27576&15250&52608&135576&55344&14118&36032&29436&115728&451584\\
\lambda_{13}  &102048&74904&137052&246192&152688&73320&116664&108072&228672&589824\\
\lambda_{14}  &111768&78088&153784&281328&171888&76240&129440&120144&263712&653184\\
\lambda_{15}  &13228&4752&35912&112752&35088&3694&19904&13788&91584&421632\\
\lambda_{16}  &63656&41176&97240&202368&108672&39412&76672&67728&182640&546048\\
\lambda_{17}  &460512&399936&512640&669792&555360&398688&486528&478176&660864&1042560\\
\lambda_{18}  &252768&207744&298032&434880&330480&205776&273984&263520&420240&799488\\
\lambda_{19}  &614880&534240&679680&869760&728640&532800&645120&640800&862560&1267200\\
\lambda_{20}  &201792&149952&256704&417408&285120&147648&225792&216576&402048&824832\\
\lambda_{21}  &20172&8284&47200&135144&48192&7024&28968&21756&113808&463104\\
\lambda_{22}  &8284&2150&28592&99576&25344&1392&13744&8688&78720&400896\\
\lambda_{23}  &47200&28592&78284&177432&85920&27228&58968&50748&157776&517248\\
\lambda_{24}  &135144&99576&177432&307128&199536&97872&153744&143544&289920&678528\\
\lambda_{25}  &48192&25344&85920&199536&92688&23376&62304&53568&179424&565632\\
\lambda_{26}  &7024&1392&27228&97872&23376&596&12120&7068&76224&398592\\
\lambda_{27}  &28968&13744&58968&153744&62304&12120&39088&31224&132096&489600\\
\lambda_{28} &21756&8688&50748&143544&53568&7068&31224&22656&120672&478080\\ 
\lambda_{29}  &113808&78720&157776&289920&179424&76224&132096&120672&269952&665856\\
\lambda_{30}   &463104&400896&517248&678528&565632&398592&489600&478080&665856&1051776 \\
\hline \end{array} 
$$}
\end{center}

\paragraph{Coxeter group of type $E_8$}
The set of representants of Coxeter classes that we consider is:\\
$
( 
\lambda_{1}=12,\,\lambda_{2}=123,\,\lambda_{3}=1234,\,\lambda_{4}=12345,\,\lambda_{5}=123456,\,\lambda_{6}=1234567,\\ 
\lambda_{7}=12345678,\,\lambda_{8}=1234568,\,\lambda_{9}=123457,\,\lambda_{10}=1234578,
\,\lambda_{11}=12346,\\
\lambda_{12}=123467,\,\lambda_{13}=1234678,\,\lambda_{14}=123468,\,\lambda_{15}=1235,\,\lambda_{16}=12356,\\
\lambda_{17}=123567,\,\lambda_{18}=1235678,\,\lambda_{19}=123568,\,\lambda_{20}=12357,
\,\lambda_{21}=1245,\\
\lambda_{22}=124567,\,
\lambda_{23}=1245678,\,\lambda_{24}=12457,\,\lambda_{25}=125,\,\lambda_{26}=1257,\\ 
\lambda_{27}=13,\,\lambda_{28}=134,\,
\lambda_{29}=13456,\,\lambda_{30}=134567,
\,\lambda_{31}=1345678,\\
\lambda_{32}=13467,\,\lambda_{33}=134678,\,\lambda_{34}=1356,\,\lambda_{35}=2345,\,\lambda_{36}=234567,\\
\lambda_{37}=2345678,\,\lambda_{38}=23457,\,\lambda_{39}=234578) .
$\\
We give here  the submatrix of $D'$ indexed by this Coxeter classes. 
We have decomposed this matrix into four tables.
\begin{center}
\tiny{$$
\begin{array}{|c|c|c|c|c|c|c|c|c|c|c|}\hline 
&\lambda_{1} &\lambda_{2} &\lambda_{3} &\lambda_{4} &\lambda_{5} &\lambda_{6} &\lambda_{7} &\lambda_{8} &\lambda_{9} &\lambda_{10} \\
\hline 
\lambda_{1} & 437160960&339793920&255744000&214179840&192015360&180034560&174182400&183121920&198650880&186301440\\ 
\lambda_{2} & 339793920&219699360&126612000&88485120&70796160&62082720&58060800&64477440&76479840&66998880\\   
\lambda_{3} & 255744000&126612000&41733360&17606160&9810000&6918480&5806080&7849680&12730800&8969280\\
\lambda_{4} & 214179840&88485120&17606160&3848832&1146672&530880&362880&767760&2308128&1171296\\ 
\lambda_{5} & 192015360&70796160&9810000&1146672&143460&29772&13440&80592&611292&236028\\ 
\lambda_{6}& 180034560&62082720&6918480&530880&29772&1596&240&15242&270192&93484\\   
\lambda_{7} & 174182400&58060800&5806080&362880&13440&240&1&6720&181440&60480\\  
\lambda_{8} & 183121920&64477440&7849680&767760&80592&15242&6720&44802&403548&151210\\   
\lambda_{9} & 198650880&76479840&12730800&2308128&611292&270192&181440&403548&1353996&658812\\  
\lambda_{10}& 186301440&66998880&8969280&1171296&236028&93484&60480&151210&658812&293848\\ 
\lambda_{11}& 228648960&102617280&27612240&10012704&5163984&3517500&2903040&4058736&7034268&4761696\\  
 \lambda_{12}& 205424640&82549440&16471440&4418016&1919700&1213572&967680&1455716&2934284&1824168\\ 
\lambda_{13}&  189527040&69662880&10325040&1800768&575808&320262&241920&413784&1105044&598522\\   
\lambda_{14}& 209064960&86107680&18873840&5760384&2723772&1788716&1451520&2102964&3929460&2548040\\  
\lambda_{15}& 289059840&165227040&78967200&49372800&37176768&31559040&29030400&33144576&41255520&34897344\\
\lambda_{16}& 244039680&118974240&41486880&20691600&13721040&10883304&9676800&11728152&16206312&12737136\\ 
\lambda_{17}&  212613120&89501760&21350160&7495728&4031292&2868280&2419200&3241876&5366368&3749372\\ 
\lambda_{18} & 192890880&72646560&12023040&2620992&1017696&617378&483840&758036&1687240&997262\\
\lambda_{19}& 220078080&97133760&27315840&11716608&7220784&5533552&4838400&6049880&8880696&6707120\\
\lambda_{20} & 252288000&128331360&50086080&27667872&19535400&16043784&14515200&17056944&22359768&18235608\\
\lambda_{21}&  272010240&144900000&58440480&30997632&20717568&16357344&14515200&17662080&24329280&19130976\\
\lambda_{22}&  202014720&79449120&14485680&3241824&1157268&644524&483840&830428&2034528&1140620\\
\lambda_{23} & 186301440&67004640&8979360&1183392&248676&103974&69120&161932&670552&304712\\ 
\lambda_{24}&  240215040&114829920&37707360&17456208&10890768&8315728&7257600&9104720&13278240&10052528\\
\lambda_{25} & 359009280&245329920&156637440&119022336&100806912&91517184&87091200&93984768&106502400&96609024\\ 
\lambda_{26} & 302561280&182062080&96975360&66309888&52937856&46522368&43545600&48278784&57293184&50215296\\
\lambda_{27} & 410572800&301190400&205804800&159148800&134887680&122135040&116121600&125556480&142352640&128977920\\
\lambda_{28} &316477440&189504000&92640000&55407360&39456000&32179200&29030400&34310400&44862720&36564480\\ 
\lambda_{29}&  217912320&92171520&20264400&5483088&2181720&1266912&967680&1591872&3538656&2079576\\
\lambda_{30}& 195425280&73762560&11420880&1850592&459240&204064&138240&304136&1076372&515792\\ 
\lambda_{31}&183121920&64434240&7825680&785184&101592&29822&17280&60852&420692&166672\\ 
\lambda_{32}&232657920&106981920&31490400&13092576&7648560&5640672&4838400&6269376&9679568&7059296\\
\lambda_{33}&205655040&82964160&16991040&4849920&2253600&1486784&1209600&1744056&3295576&2127960\\
\lambda_{34}&276480000&150465600&64491840&36772128&26098848&21407328&19353600&22766976&29768256&24300000\\ 
\lambda_{35}&251596800&121875840&37313280&14054400&6943104&4503168&3628800&5326848&9694080&6350976\\ 
\lambda_{36}&192199680&71055360&10046880&1273920&198048&55904&30240&117144&694728&284088\\
\lambda_{37}&180126720&62232480&7044960&581664&43704&5694&2160&23510&300392&107758\\ 
\lambda_{38}&225054720&99004320&24805440&8013312&3656160&2289504&1814400&2759136&5386368&3393840\\
\lambda_{39}& 202245120&79833600&14928480&3561984&1361136&789948&604800&993696&2276868&1320864\\ 
\hline \end{array} 
$$}
\end{center}
\begin{center}
\tiny{$$
\begin{array}{|c|c|c|c|c|c|c|c|c|c|c|}\hline 
&\lambda_{11} &\lambda_{12} &\lambda_{13} &\lambda_{14} &\lambda_{15} &\lambda_{16} &\lambda_{17} &\lambda_{18} &\lambda_{19} &\lambda_{20} \\
\hline 
\lambda_{1}&  228648960&205424640&189527040&209064960&289059840&244039680&212613120&192890880&220078080&252288000\\
\lambda_{2} & 102617280&82549440&69662880&86107680&165227040&118974240&89501760&72646560&97133760&128331360\\ 
\lambda_{3} & 27612240&16471440&10325040&18873840&78967200&41486880&21350160&12023040&27315840&50086080\\
\lambda_{4} & 10012704&4418016&1800768&5760384&49372800&20691600&7495728&2620992&11716608&27667872\\
\lambda_{5}&  5163984&1919700&575808&2723772&37176768&13721040&4031292&1017696&7220784&19535400\\ 
\lambda_{6} &3517500&1213572&320262&1788716&31559040&10883304&2868280&617378&5533552&16043784\\ 
\lambda_{7} & 2903040&967680&241920&1451520&29030400&9676800&2419200&483840&4838400&14515200\\ 
\lambda_{8} & 4058736&1455716&413784&2102964&33144576&11728152&3241876&758036&6049880&17056944\\          
\lambda_{9} & 7034268&2934284&1105044&3929460&41255520&16206312&5366368&1687240&8880696&22359768\\ 
\lambda_{10}&4761696&1824168&598522&2548040&34897344&12737136&3749372&997262&6707120&18235608\\ 
\lambda_{11} & 17498496&9734340&5677848&11465916&60756288&29139360&13445220&6837720&18242544&36592632\\          
 \lambda_{12}&  9734340&4762868&2383000&5924412&45895296&19414512&7478296&3101068&11268400&25837632\\          
\lambda_{13} & 5677848&2383000&948816&3182820&36815232&13934304&4443900&1404342&7543112&19588272\\          
\lambda_{14} &11465916&5924412&3182820&7197180&48708288&21425424&8819584&3992964&12776752&27992592\\          
\lambda_{15} &60756288&45895296&36815232&48708288&115584864&75212064&51515040&39026592&57994176&83879424\\  
\lambda_{16}& 29139360&19414512&13934304&21425424&75212064&41410752&23615184&15378384&28820496&49126992\\ 
\lambda_{17} & 13445220&7478296&4443900&8819584&51515040&23615184&10498956&5317000&14573384&30287568\\          
\lambda_{18}  & 6837720&3101068&1404342&3992964&39026592&15378384&5317000&1930534&8569664&21186744\\          
\lambda_{19}&  18242544&11268400&7543112&12776752&57994176&28820496&14573384&8569664&18901640&35693808\\          
\lambda_{20} & 36592632&25837632&19588272&27992592&83879424&49126992&30287568&21186744&35693808&56943600\\ 
\lambda_{21}& 41920320&28605120&20803008&31287360&95738400&56280480&33963456&22816416&40300608&65026656\\           
\lambda_{22} & 8273472&3741876&1643036&4821348&43491840&17712384&6318552&2293132&9973912&24004224\\          
\lambda_{23} & 4773192&1836244&609716&2559460&34903776&12747216&3760768&1008216&6716768&18242856\\          
\lambda_{24} & 25772256&16401888&11196960&18379008&71406624&37965408&20531616&12598624&25688512&45662208\\   
\lambda_{25} & 132655104&112520448&99330048&115941888&192948480&148253184&119280384&102279168&126604800&157102848\\ 
\lambda_{26} &77854848&62131968&52286976&64998528&132725760&92219520&67873152&54608640&74390016&100761600\\           
\lambda_{27}  &175737600&149921280&132503040&154172160&244166400&193363200&158112000&136339200&166717440&202901760\\          
\lambda_{28}&  69719040&50734080&39014400&54366720&133651200&86676480&57722880&41913600&65579520&96606720\\ 
\lambda_{29} & 11999184&5813304&2791248&7258320&52382496&22940424&9063408&3704664&13435968&30030528\\          
\lambda_{30}&  6229488&2529564&915192&3432572&39367392&15106488&4810340&1432848&8168688&21089976\\          
\lambda_{31}&4057860&1468100&428652&2111844&33130848&11729184&3251596&771600&6059208&17058936\\          
\lambda_{32} &20798256&12513072&8054880&14307584&64718688&32587152&16345568&9287520&21236736&40084896\\          
\lambda_{33}&10187232&5144832&2701152&6313696&46364160&19874784&7876864&3430720&11676352&26291712\\
\lambda_{34}&47684160&34102080&26024832&36765504&101408256&61969248&39464928&28053504&45759168&70604352\\
\lambda_{35}& 23837184&13286016&7615872&15634944&74620800&37613952&18058752&9239040&23968512&46244736\\
\lambda_{36}&5339088&2029112&636624&2847288&37450368&13938816&4168160&1092280&7383336&19768752\\ 
\lambda_{37}& 3601860&1258932&341066&1843456&31700928&10988856&2929372&645916&5608648&16158480\\          
\lambda_{38}&15251184&7964928&4245936&9629280&57751776&26730960&11561520&5340288&16268256&34126176\\ 
\lambda_{39}& 8629560&4008432&1838640&5103384&43913952&18091512&6611700&2505792&10290600&24394272\\ 
\hline \end{array} 
$$}
\end{center}
\begin{center}
\tiny{$$
\begin{array}{|c|c|c|c|c|c|c|c|c|c|c|}\hline 
&\lambda_{21} &\lambda_{22} &\lambda_{23} &\lambda_{24} &\lambda_{25} &\lambda_{26} &\lambda_{27} &\lambda_{28} &\lambda_{29} &\lambda_{30} \\
\hline 
\lambda_{1}&   272010240&202014720&186301440&240215040&359009280&302561280&410572800&316477440&217912320&195425280\\          
\lambda_{2} & 144900000&79449120&67004640&114829920&245329920&182062080&301190400&189504000&92171520&73762560\\
\lambda_{3} & 58440480&14485680&8979360&37707360&156637440&96975360&205804800&92640000&20264400&11420880\\
\lambda_{4} & 30997632&3241824&1183392&17456208&119022336&66309888&159148800&55407360&5483088&1850592\\
\lambda_{5}&  20717568&1157268&248676&10890768&100806912&52937856&134887680&39456000&2181720&459240\\ 
\lambda_{6}& 16357344&644524&103974&8315728&91517184&46522368&122135040&32179200&1266912&204064\\ 
\lambda_{7}&  14515200&483840&69120&7257600&87091200&43545600&116121600&29030400&967680&138240\\
\lambda_{8}&  17662080&830428&161932&9104720&93984768&48278784&125556480&34310400&1591872&304136\\          
\lambda_{9} & 24329280&2034528&670552&13278240&106502400&57293184&142352640&44862720&3538656&1076372\\   
\lambda_{10}&19130976& 1140620&304712&10052528&96609024&50215296&128977920&36564480&2079576&515792\\          
\lambda_{11} &41920320&8273472&4773192&25772256&132655104&77854848&175737600&69719040&11999184&6229488\\          
 \lambda_{12}& 28605120&3741876&1836244&16401888&112520448&62131968&149921280&50734080&5813304&2529564\\          
\lambda_{13}&  20803008&1643036&609716&11196960&99330048&52286976&132503040&39014400&2791248&915192\\          
\lambda_{14}& 31287360&4821348&2559460&18379008&115941888&64998528&154172160&54366720&7258320&3432572\\
\lambda_{15}& 95738400&43491840&34903776&71406624&192948480&132725760&244166400&133651200&52382496&39367392\\          
\lambda_{16}& 56280480&17712384&12747216&37965408&148253184&92219520&193363200&86676480&22940424&15106488\\ 
\lambda_{17} & 33963456&6318552&3760768&20531616&119280384&67873152&158112000&57722880&9063408&4810340\\          
\lambda_{18}&  22816416&2293132&1008216&12598624&102279168&54608640&136339200&41913600&3704664&1432848\\          
\lambda_{19}&  40300608&9973912&6716768&25688512&126604800&74390016&166717440&65579520&13435968&8168688\\          
\lambda_{20} & 65026656&24004224&18242856&45662208&157102848&100761600&202901760&96606720&30030528&21089976\\          
\lambda_{21} & 75382176&26365632&19138656&52452096&173882880&113431680&224570880&111989760&33873792&22662144\\          
\lambda_{22} & 26365632&2781016&1153396&14736768&109458432&59635968&146085120&47696640&4562088&1696248\\          
\lambda_{23} & 19138656&1153396&316042&10062240&96609024&50216832&128977920&36564480&2094384&529212\\   
\lambda_{24} &  52452096&14736768&10062240&34527840&144317952&88473216&189008640&82394880&19656384&12217728\\
\lambda_{25}& 173882880&109458432&96609024&144317952&269266176&208760832&323066880&217232640&122595840&103799808\\            
\lambda_{26}&  113431680&59635968&50216832&88473216&208760832&149380608&259614720&151741440&69357312&55273728\\ 
\lambda_{27}& 224570880&146085120&128977920&189008640&323066880&259614720&380816640&274648320&163296000&138620160\\           
\lambda_{28}&  111989760&47696640&36564480&82394880&217232640&151741440&274648320&157002240&59082240&42251520\\          
\lambda_{29}&   33873792&4562088&2094384&19656384&122595840&69357312&163296000&59082240&7227744&2992704\\          
\lambda_{30}& 22662144&1696248&529212&12217728&103799808&55273728&138620160&42251520&2992704&858964\\          
\lambda_{31}&17628864&845324&177668&9098832&93987072&48285696&125452800&34218240&1604760&322976\\ 
\lambda_{32}&45871008&10992384&7072288&29190720&136832256&81780096&180195840&74215680&15152928&8804144\\          
\lambda_{33}&29095296&4115744&2140864&16857120&112909824&62584320&150128640&51194880&6258240&2881408\\          
\lambda_{34}&81255744&31848480&24312672&58146048&179071488&118874880&229547520&117849600&39694752&28099584\\          
\lambda_{35}&54028800&11355264&6352128&33866496&152239104&92799360&201139200&87736320&16608384&8449920\\
\lambda_{36}&20966976&1254496&297368&11094080&101090304&53230848&134991360&39705600&2326656&531624\\          
\lambda_{37}& 16495872&682704&118388&8417920&91657728&46665984&122238720&32348160&1329144&228012\\ 
\lambda_{38}& 38966976&6556128&3402288&23387376&129238272&74884992&171590400&66055680&9916080&4642992\\          
\lambda_{39}& 26796672&3034584&1332756&15110736&109840896&60064128&146292480&48142080&4887600&1918056\\
\hline \end{array} 
$$}
\end{center}

\begin{center}
\tiny{$$
\begin{array}{|c|c|c|c|c|c|c|c|c|c|}\hline 
&\lambda_{31} &\lambda_{32} &\lambda_{33} &\lambda_{34} &\lambda_{35} &\lambda_{36} &\lambda_{37} &\lambda_{38} &\lambda_{39} \\
\hline 
\lambda_{1} & 183121920&232657920&205655040&276480000&251596800&192199680&180126720&225054720&202245120\\
\lambda_{2} & 64434240&106981920&82964160&150465600&121875840&71055360&62232480&99004320&79833600\\
\lambda_{3} &  7825680&31490400&16991040&64491840&37313280&10046880&7044960&24805440&14928480\\
\lambda_{4} & 785184&13092576&4849920&36772128&14054400&1273920&581664&8013312&3561984\\
\lambda_{5} & 101592&7648560&2253600&26098848&6943104&198048&43704&3656160&1361136\\
\lambda_{6}& 29822&5640672&1486784&21407328&4503168&55904&5694&2289504&789948\\
\lambda_{7} & 17280&4838400&1209600&19353600&3628800&30240&2160&1814400&604800\\
\lambda_{8} &  60852&6269376&1744056&22766976&5326848&117144&23510&2759136&993696\\
\lambda_{9} & 420692&9679568&3295576&29768256&9694080&694728&300392&5386368&2276868\\
\lambda_{10}&166672& 7059296&2127960&24300000&6350976&284088&107758&3393840&1320864\\
\lambda_{11}& 4057860&20798256&10187232&47684160&23837184&5339088&3601860&15251184&8629560\\
 \lambda_{12}& 1468100&12513072&5144832&34102080&13286016&2029112&1258932&7964928&4008432\\
\lambda_{13}&  428652&8054880&2701152&26024832&7615872&636624&341066&4245936&1838640\\
\lambda_{14}& 2111844&14307584&6313696&36765504&15634944&2847288&1843456&9629280&5103384\\
\lambda_{15} &33130848&64718688&46364160&101408256&74620800&37450368&31700928&57751776&43913952\\
\lambda_{16}& 11729184&32587152&19874784&61969248&37613952&13938816&10988856&26730960&18091512\\
\lambda_{17}&  3251596&16345568&7876864&39464928&18058752&4168160&2929372&11561520&6611700\\
\lambda_{18}&  771600&9287520&3430720&28053504&9239040&1092280&645916&5340288&2505792\\
\lambda_{19} &6059208&21236736&11676352&45759168&23968512&7383336&5608648&16268256&10290600\\
\lambda_{20}& 17058936&40084896&26291712&70604352&46244736&19768752&16158480&34126176&24394272\\ 
\lambda_{21} &  17628864&45871008&29095296&81255744&54028800&20966976&16495872&38966976&26796672\\
\lambda_{22} & 845324&10992384&4115744&31848480&11355264&1254496&682704&6556128&3034584\\
\lambda_{23} & 177668&7072288&2140864&24312672&6352128&297368&118388&3402288&1332756\\
\lambda_{24} &  9098832&29190720&16857120&58146048&33866496&11094080&8417920&23387376&15110736\\
\lambda_{25}&  93987072&136832256&112909824&179071488&152239104&101090304&91657728&129238272&109840896\\
\lambda_{26} &  48285696&81780096&62584320&118874880&92799360&53230848&46665984&74884992&60064128\\
\lambda_{27} &125452800&180195840&150128640&229547520&201139200&134991360&122238720&171590400&146292480\\ 
\lambda_{28}& 34218240&74215680&51194880&117849600&87736320&39705600&32348160&66055680&48142080\\
\lambda_{29}&   1604760&15152928&6258240&39694752&16608384&2326656&1329144&9916080&4887600\\
\lambda_{30}& 322976&8804144&2881408&28099584&8449920&531624&228012&4642992&1918056\\
\lambda_{31}&76360&6271344&1757552&22752288&5300352&138008&38334&2758608&1007208\\
\lambda_{32}&6271344&24146240&12973696&51653952&27668736&7839456&5731696&18505056&11352048\\
\lambda_{33}&1757552&12973696&5530880&34602912&13800960&2368040&1534416&8413056&4383360\\
\lambda_{34}& 22752288&51653952&34602912&87158304&60018048&26369280&21546720&44713152&32272992\\
\lambda_{35}&5300352&27668736&13800960&60018048&33041664&7168896&4625280&21125376&11822976\\
\lambda_{36}&138008&7839456&2368040&26369280&7168896&256784&72328&3819456&1460880\\
\lambda_{37}&38334&5731696&1534416&21546720&4625280&72328&11002&2368224&830520\\
\lambda_{38}&2758608&18505056&8413056&44713152&21125376&3819456&2368224&13069968&6915312\\
\lambda_{39}& 1007208&11352048&4383360&32272992&11822976&1460880&830520&6915312&3292116\\
\hline \end{array} 
$$}
\end{center}

\subsection*{Acknowledgments}

The author is indebted to C.~Reutenauer for his invaluable
contribution during this work. He also thanks M.~Schocker for the many
hours passed on the double coset conjecture; and R.B~Howlett for our
useful discussions. Finally, he is grateful to F.~L{\"u}beck for his
helpful remarks about CHEVIE and GAP.


\begin{thebibliography}{99}

\bibitem{bergeron} F.~Bergeron, N.~Bergeron, R.~B.~Howlett and D.~E.~Taylor: \emph{A decomposition of the descent algebra of a finite Coxeter group.},
J. Algebraic Combin., 1992, \textbf{1}, p.23-44.
\bibitem{nantel} N.~Bergeron: \emph{A Decomposition of the Descent Algebra of the Hyperoctahedral Group II}, J. Algebra, 1992, \textbf{148}, p.98-122.
\bibitem{bourbaki} N.~Bourbaki: \emph{Groupes et alg\`ebres de Lie}, Chap. 4-6, Hermann, 1968.
\bibitem{fleisch} P.~Fleischmann: \emph{On pointwise conjugacy of distinguished coset representatives in Coxeter groups}, J. Group Theory, 2002, \textbf{5}, 
p.269 -283.
\bibitem{fulton} W.~Fulton, J.~Harris: \emph{Representation theory, a first course}, Springer, 1999, \textbf{129}. 
\bibitem{garsia} A.~M.~Garsia, C.~Reutenauer: \emph{A decomposition of Solomon's descent algebra}, Adv. in Math., 1989, \textbf{77}, p.189-262.
\bibitem{chevie} M.~Geck, G.~Hiss, F.~L{\"u}beck, G.~Malle, and G.~Pfeiffer.
\newblock {\sf CHEVIE} -- {A} system for computing and processing generic
  character tables for finite groups of {L}ie type, {W}eyl groups and {H}ecke
  algebras.
\newblock {\em Appl. Algebra Engrg. Comm. Comput.}, 7:175--210, 1996.        
\bibitem{geck} M.~Geck, G.~Pfeiffer: \emph{Characters of Finite Coxeter Groups and Iwahori-Hecke Algebras}, London Math. Soc. Mon. New Series, 2000,
 \textbf{21}.
 \bibitem{gessel} I.~Gessel: \emph{Multipartite P-partitions and inner product of skew Schur functions}, 
1984, Contemp. Math., \textbf{34}, p.289-301.
\bibitem{gesselreut} I.~Gessel, C.~Reutenauer: \emph{Counting permutations with given cycle structure and descent set},
 1993, J. Combin. Thery Ser. A, \textbf{64}, p.189-215.
\bibitem{hohlweg} C.~Hohlweg, M.~Schocker: \emph{On Deodhar's Lemma and a parabolic symmetry of finite Coxeter groups}, Preprint Institut de Recherche 
Math\'ematique avanc\'ee (IRMA), 2002.
\bibitem{humphreys} J.E.~Humphreys: \emph{Reflection groups and Coxeter groups}, 1990, Cambridge university press, \textbf{29}.
\bibitem{jollen} A.~J\" ollenbeck, C.~Reutenauer: \emph{Eine Symmetrieeigenschaft von Solomons 
Algebra und der h\" oheren Lie-Charaktere}, Abh. Math. Sem. Univ. Hamburg, \textbf{71}, 2001, P.105-111.
\bibitem{ledermann} W.~Ledermann: \emph{Introduction to group characters}, 1977, Cambridge university press.
\bibitem{poirier} S.~Poirier: \emph{Cycle type and descent set in wreath products},
 1998, Discrete Math., \textbf{180}, p.315-343.
\bibitem{reutenau} C.~Reutenauer: \emph{Free Lie algebras}, 1993,
London Math. Soc. Mon. New Series, \textbf{7}.
\bibitem{gap}
              Martin Sch{\accent127 o}nert et~al.
    \newblock {\em {GAP} -- {Groups}, {Algorithms}, and {Programming}}.
    \newblock Lehrstuhl D f{\accent127 u}r Mathematik,
              Rheinisch Westf{\accent127 a}lische Technische Hoch\-schule,
              Aachen, Germany, fifth edition, 1995. 
\bibitem{solomon68} L.~Solomon: \emph{A Decomposition of the Group Algebra of a finite Coxeter group}, 1968, J. Algebra, \textbf{9}, No $2$, p.220-239.           
\bibitem{solomon} L.~Solomon: \emph{A Mackey formula in the group ring of a Coxeter group}, 1976, J. Algebra, \textbf{41}, p.255-268.


\end{thebibliography}
\end{document}